\documentclass[11pt]{article}
\usepackage[top=1in, bottom=1in, left=1in, right=1in]{geometry}
\usepackage{tikz}
\usetikzlibrary{shapes,arrows,matrix,patterns,positioning}
\usepackage{color}
\usepackage{graphicx}
\usepackage{subcaption}
\usepackage{algorithmic}
\usepackage{amssymb}
\usepackage{epstopdf}
\usepackage[ruled]{algorithm2e}
\usepackage{float}
\usepackage{amsmath,amsthm,bm,color,epsfig,enumerate,caption}
\usepackage{hyperref}
\usepackage{booktabs}
\usepackage{multirow}
\usepackage{array}
\usepackage{bigstrut}
\usepackage{mathabx}

\newcommand{\norm}[1]{\left\lVert#1\right\rVert}
\renewcommand{\O}{O}

\newcommand{\bbR}{\mathbb{R}}
\newcommand{\calC}{\mathcal{C}}
\newcommand{\calI}{\mathcal{I}}
\newcommand{\calF}{\mathcal{F}}

\renewcommand{\vec}[1]{{\mathbf #1}}
\newcommand{\vecx}{\vec{x}}
\newcommand{\vecj}{\vec{j}}
\newcommand{\vece}{\vec{e}}
\newcommand{\vecZero}{\vec{0}}
\newcommand{\hF}{\widehat{F}}
\newcommand{\cF}{\widebar{\widebar{F}}}
\newcommand{\wtA}{\widetilde{A}}
\newcommand{\wtB}{\widetilde{B}}
\newcommand{\wtD}{\widetilde{D}}
\newcommand{\wts}{\widetilde{\Sigma}}
\newcommand{\wideT}[1]{{\left( #1 \right)^T}}
\newcommand{\wideInv}[1]{{\left( #1 \right)^{-1}}}
\newcommand{\wideInvT}[1]{{\left( #1 \right)^{-T}}}
\newcommand{\floor}[1]{{\lfloor #1 \rfloor}}
\newcommand{\defeq}{\equiv}

\makeatletter
\newcommand*{\extendadd}{
  \mathbin{
    \mathpalette\extend@add{}
  }
}
\newcommand*{\extend@add}[2]{
  \ooalign{
    $\m@th#1\leftrightarrow$%
    \vphantom{$\m@th#1\updownarrow$}
    \cr
    \hfil$\m@th#1\updownarrow$\hfil
  }
}
\makeatother

\begin{document}

\title{Distributed-memory Hierarchical Interpolative Factorization}

\author{Yingzhou Li$^\sharp$,
    Lexing Ying$^{\dagger\sharp}$
  \vspace{0.1in}\\
  $\sharp$ ICME, Stanford University\\
  $\dagger$ Department of Mathematics, Stanford University\\
}

\maketitle

\begin{abstract}
    The hierarchical interpolative factorization ({HIF}) offers an efficient
    way for solving or preconditioning elliptic partial differential
    equations. By exploiting locality and low-rank properties of the
    operators, the HIF achieves quasi-linear complexity for factorizing the
    discrete positive definite elliptic operator and linear complexity for
    solving the associated linear system.  In this paper, the
    distributed-memory {HIF} ({DHIF}) is introduced as a parallel and
    distributed-memory implementation of the HIF. The DHIF organizes the
    processes in a hierarchical structure and keep the communication as
    local as possible.  The computation complexity is $\O\left(\frac{N\log
    N}{P}\right)$ and $\O\left(\frac{N}{P}\right)$ for constructing and
    applying the DHIF, respectively, where $N$ is the size of the problem
    and $P$ is the number of processes.  The communication complexity is
    $\O\left(\sqrt{P}\log^3 P\right)\alpha +
    \O\left(\frac{N^{2/3}}{\sqrt{P}}\right)\beta$ where $\alpha$ is the
    latency and $\beta$ is the inverse bandwidth. Extensive numerical
    examples are performed on the {NERSC Edison} system with up to 8192
    processes. The numerical results agree with the complexity analysis and
    demonstrate the efficiency and scalability of the {DHIF}.
\end{abstract}



\section{Introduction}
\label{sec:Introduction}

This paper proposes an efficient distributed-memory algorithm for solving
elliptic partial differential equations (PDEs) of the form,
\begin{equation}\label{eq:original-prob}
    -\nabla \cdot (a(x)\nabla u(x)) + b(x)u(x) = f(x)
    \quad x\in \Omega\subset \bbR^3,
\end{equation}
with a certain boundary condition, where $a(x)>0$, $b(x)$ and $f(x)$ are
given functions and $u(x)$ is an unknown function. Since this elliptic
equation is of fundamental importance to problems in physical sciences,
solving \eqref{eq:original-prob} effectively has a significant impact in
practice. Discretizing this with local schemes such as the finite difference
or finite element methods leads to a sparse linear system,
\begin{equation}\label{eq:discretized-prob}
    Au = f,
\end{equation}
where $A\in \bbR^{N\times N}$ is a sparse symmetric matrix with $\O(N)$
non-zero entries with $N$ being the number of the discretization points, and
$u$ and $f$ are the discrete approximations of the functions $u(x)$ and
$f(x)$, respectively. For many practical applications, one often needs to
solve \eqref{eq:original-prob} on a sufficient fine mesh for which $N$ can
be very large, especially for three dimensional (3D) problems. Hence, there
is a practical need for developing fast and parallel algorithms for the
efficient solution of \eqref{eq:original-prob}.

\subsection{Previous work}
\label{sub:Previous work}

A great deal of effort in the field of scientific computing has been devoted
to the efficient solution of \eqref{eq:discretized-prob}.  Beyond the
$\O(N^3)$ complexity na\"{i}ve matrix inversion approach, one can classify
the existing fast algorithms into the following groups.

The first one consists of the sparse direct algorithms, which take advantage
of the sparsity of the discrete problem.  The most noticeable example in
this group is the nested dissection multifrontal method (MF)
method~\cite{George1973,Duff1983,Liu1992}. By carefully exploring the
sparsity and the locality of the problem, the multifrontal method factorizes
the matrix $A$ (and thus $A^{-1}$) as a product of sparse lower and upper
triangular matrices. For 3D problems, the factorization step costs $\O(N^2)$
operations while the application step takes $\O(N^{4/3})$ operations. Many
parallel implementations~\cite{Amestoy2000,Amestoy2001,Poulson2013a} of the
multifrontal method were proposed and they typically work quite well for
problem of moderate size. However, as the problem size goes beyond a couple
of millions, most implementations, including the distributed-memory ones,
hit severe bottlenecks in memory consumption.

The second group consists of iterative
solvers~\cite{Briggs2000,Saad2001,Saad2003,Falgout2000}, including famous
algorithms such as the conjugate gradient (CG) method and the multigrid
method. Each iteration of these algorithms typically takes $O(N)$ steps and
hence the overall cost for solving \eqref{eq:discretized-prob} is
proportional to the number of iterations required for convergence. For
problems with smooth coefficient functions $a(x)$ and $b(x)$, the number of
iterations typically remains rather small and the optimal linear complexity
is achieved.  However, if the coefficient functions lack regularity or have
high contrast, the iteration number typically grows quite rapidly as the
problem size increases.

The third group contains the methods based on structured matrices
\cite{Bebendorf2003,Bebendorf2005,Borm2010,Chandrasekaran2010}. These
methods, for example the
$\mathcal{H}$-matrix~\cite{Hackbusch1999,Hackbusch2000}, the
$\mathcal{H}^2$-matrix~\cite{Hackbusch2002}, and the hierarchically
semi-separable matrix (HSS)~\cite{Chandrasekanran2005,Xia2010}, are shown to
have efficient approximations of linear or quasi-linear complexity for the
matrices $A$ and $A^{-1}$. As a result, the algebraic operations of these
matrices are of linear or quasi-linear complexities as well. More
specifically, the recursive inversion and the rocket-style
inversion~\cite{Ambikasaran2013} are two popular methods for the inverse
operation. For distributed-memory implementations, however, the former lacks
parallel scalability~\cite{Izadi2012,Kriemann2013} while the latter
demonstrates scalability only for 1D and 2D problems~\cite{Ambikasaran2013}.
For 3D problems, these methods typically suffer from large prefactors that
make them less efficient for practical large-scale problems.

A recent group of methods explore the idea of integrating the MF method with
the hierarchical matrix
\cite{Martinsson2009,Xia2009,Xia2013,Xia2013a,Gillman2014,Wang2016,Hao2016}
or block low-rank matrix~\cite{Schmitz2012,Schmitz2014,Amestoy2015} approach
in order to leverage the efficiency of both methods. Instead of directly
applying the hierarchical matrix structure to the 3D problems, these methods
apply it to the representation of the frontal matrices (i.e., the
interactions between the lower dimensional fronts). These methods are of
linear or quasi-linear complexities in theory with much small prefactors.
However, due to the combined complexity, the implementation is highly
non-trivial and quite difficult for parallelization~\cite{Liu2016,Xin2014}.

More recently, the hierarchical interpolative factorization
(HIF)~\cite{Ho2015,Ho2016} is proposed as a new way for solving elliptic
PDEs and integral equations. As compared to the multifrontal method, the HIF
includes an extra step of skeletonizing the fronts in order to reduce the
size of the dense frontal matrices. Based on the key observation that the
number of skeleton points on each front scales linearly as the
one-dimensional fronts, the HIF factorizes the matrix $A$ (and thus
$A^{-1}$) as a product of sparse matrices that contains only $\O(N)$
non-zero entries in total. In addition, the factorization and application of
the HIF are of complexities $\O(N\log N)$ and $\O(N)$, respectively, for $N$
being the total number of degrees of freedom (DOFs) in
\eqref{eq:discretized-prob}. In practice, the HIF shows significant saving
in terms of computational resources required for 3D problems.

\subsection{Contribution}
\label{sub:Contribution}

This paper proposes the first {\em distributed-memory} hierarchical
interpolative factorization (DHIF) for solving very large scale problems. In
a nutshell, the DHIF organizes the processes in an octree structure in the
same way that the HIF partitions the computation domain. In the simplest
setting, each leaf node of the computation domain is assigned a single
process. Thanks to the locality of the operator in \eqref{eq:original-prob},
this process only communicates with its neighbors and all algebraic
computations are local within the leaf node. At higher levels, each node of
the computation domain is associated with a process group that contains all
processes in the subtree starting from this node. The computations are all
local within this process group via parallel dense linear algebra and the
communications are carried out between neighboring process groups. By
following this octree structure, we make sure that both the communication
and computations in the distributed-memory HIF are evenly distributed. As a
result, the distributed-memory HIF implementation achieves
$\O\left(\frac{N\log N}{P}\right)$ and $\O\left(\frac{N}{P}\right)$ parallel
complexity for constructing and applying the factorization, respectively,
where $N$ is the number of DOFs and $P$ is the number of processes.

We have performed extensive numerical tests. The numerical results support
the complexity analysis of the distributed-memory HIF and suggest that the
DHIF is a scalable method up to thousands of processes and can be applied to
solve large scale elliptic PDEs.

\subsection{Organization}
\label{sub:Organization}

The rest of this paper is organized as follow.  In
Section~\ref{sec:Preliminaries}, we review the basic tools needed for both
HIF and DHIF, and review the sequential HIF.
Section~\ref{sec:Distributed-memory HIF} presents the DHIF as a parallel
extension of the sequential HIF for 3D problems. Complexity analyses for
memory usage, computation time and communication volume are given at the end
of this section.  The numerical results detailed in
Section~\ref{sec:Numerical Results} show that the DHIF is applicable to
large scale problems and achieves parallel scalability up to thousands of
processes. Finally, Section~\ref{sec:Conclusion} concludes with some extra
discussions on future work.

\section{Preliminaries}
\label{sec:Preliminaries}

This section reviews the basic tools and the sequential HIF. First, we start
by listing the notations that are widely used throughout this paper.

\subsection{Notations}
\label{sub:Notations}

In this paper, we adopt MATLAB notations for simple representation of
submatrices.  For example, given a matrix $A$ and two index sets, $s_1$ and
$s_2$, $A(s_1,s_2)$ represents the submatrix of $A$ with the row indices in
$s_1$ and column indices in $s_2$.  The next two examples explores the usage
of MATLAB notation ``$:$".  With the same settings, $A(s_1,:)$ represents
the submatrix of $A$ with row indices in $s_1$ and all columns.  Another
usage of notation ``$:$" is to create regularly spaced vectors for integer
values $i$ and $j$, for instance, $i:j$ is the same as $[i,i+1,i+2,\dots,j]$
for $i\leq j$.

In order to simplify the presentation, we consider the problem
\eqref{eq:original-prob} with periodic boundary condition and assume that
the domain $\Omega = [0,1)^3$, and is discretized with a grid of size
$n\times n\times n$ for $n=2^Lm$, where $L=\O(\log n)$ and $m=\O(1)$ are
both integers. In the rest of this paper, $L+1$ is known as the number of
levels in the hierarchical structure and $L$ is the level number of the root
level.  We use $N = n^3$ to denote the total number of DOFs, which is the
dimension of the sparse matrix $A$ in \eqref{eq:discretized-prob}.
Furthermore, each grid point $\vecx_\vecj$ is defined as
\begin{equation}
  \vecx_\vecj = h\vecj = h(j_1,j_2,j_3)
\end{equation}
where $h = 1/n$, $\vecj = (j_1,j_2,j_3)$ and $0\leq j_1,j_2,j_3 < n$.

\begin{figure}[htp]
    \centering
    \includegraphics[height=1.5in]{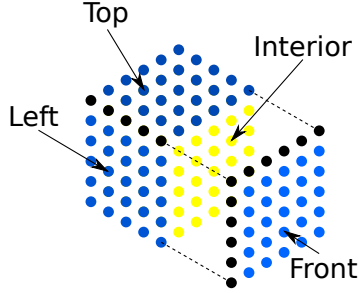}
    \caption{Cell Structure: top, front, left, and interior points are
      indicated by arrows; bottom, back, and right points are not
      plotted in the figure; the black dots denote the edge points;
      the dash line indicates that the front frame is pulled away in
      order to show the interior points.}
    \label{fig:cellstructure}
\end{figure}

In order to fully explore the hierarchical structure of the problem, we
recursively bipartite each dimension of the grid into $L+1$ levels.  Let the
leaf level be level $0$ and the root level be level $L$.  At level $\ell$, a
cell indexed with $\vecj$ is of size $m2^\ell \times m2^\ell \times m2^\ell$
and each point in the cell is in the range, $\left[m2^\ell j_1+(0:m2^\ell
-1) \right] \times \left[m2^\ell
  j_2+(0:m2^\ell -1) \right] \times \left[m2^\ell j_3+(0:m2^\ell -1)
  \right],$ for $\vecj = (j_1,j_2,j_3)$ and $0\leq
j_1,j_2,j_3<2^{L-\ell}$.  $C^\ell_\vecj$ denotes the grid point set of the
cell at level $\ell$ indexed with $\vecj$.

A cell $C^\ell_\vecj$ owns three faces: top, front, and left. Each of these
three faces contains the grid points on the first frame in the corresponding
direction.  For example, the front face contains the grid points in
$\left[m2^\ell j_1+(0:m2^\ell -1) \right] \times \left[m2^\ell j_2 \right]
\times \left[m2^\ell j_3+(0:m2^\ell -1)
  \right]$.  Besides these three in-cell faces (top, front, and left)
that are owned by a cell, each cell is also adjacent to three out-of-cell
faces (bottom, back, right) owned by its neighbors. Each of these three
faces contains the grid points on the next to the last frame in the
corresponding dimension. As a result, these faces contain DOFs that belong
to adjacent cells. For example, the bottom face of $C^\ell_\vecj$ contains
the grid points in $\left[m2^\ell (j_1+1)
  \right] \times \left[m2^\ell j_2+(0:m2^\ell -1) \right] \times
\left[m2^\ell j_3+(0:m2^\ell -1) \right]$. These six faces are the
surrounding faces of $C^\ell_\vecj$.  One also defines the interior of
$C^\ell_\vecj$ to be $I^\ell_{\vecj} = \left[m2^\ell j_1+(1:m2^\ell
  -1) \right] \times \left[m2^\ell j_2+(1:m2^\ell -1) \right] \times
\left[m2^\ell j_3+(1:m2^\ell -1) \right]$ for the same $\vecj =
(j_1,j_2,j_3)$ and $0\leq j_1,j_2,j_3 < 2^{L-\ell}$.  Figure
\ref{fig:cellstructure} gives an illustration of a cell, its faces, and its
interior. These definitions and notations are summarized in
Table~\ref{tab:Notations}. Also included here are some notations used for
the processes, which will be introduced later.

\begin{table}[ht!]
  \centering
  \begin{tabular}{l|l}
    \toprule
    Notation & Description\\
    \toprule
    $n$ &    Number of points in each dimension of the grid\\
    $N$ &    Number of points in the grid\\
    $h$ &    Grid gap size\\
    $\ell$ &    Level number in the hierarchical structure\\
    $L$ &    Level number of the root level in the hierarchical structure\\
    \midrule
    $\vece_1$, $\vece_2$, $\vece_3$ &    Unit vector along each dimension\\
    $\vecZero$ &    Zero vector\\
    $\vecj$ &    Triplet index $\vecj=(j_1,j_2,j_3)$\\
    $\vecx_\vecj$ &    Point on the grid indexed with $\vecj$\\
    \midrule
    $\Omega$ &    The set of all points on the grid\\
    $C^\ell_\vecj$ &    Cell at level $\ell$ with index $\vecj$\\
    $\calC^\ell$ &    $\calC^\ell = \{C^\ell_\vecj\}_\vecj$ is the set of all cells at level $\ell$\\
    $\calF^\ell_\vecj$ &    Set of all surrounding faces of cell $C^\ell_\vecj$\\
    $\calF^\ell$ &    Set of all faces at level $\ell$\\
    $I^\ell_\vecj$ &    Interior of $C^\ell_\vecj$\\
    $\calI^\ell$ &    $\calI^\ell = \{I^\ell_\vecj\}_\vecj$ is the set of all interiors    at level $\ell$\\
    \midrule
    $\Sigma^\ell$ &    The set of active DOFs at level $\ell$\\
    $\Sigma^\ell_\vecj$ &    The set of active DOFs at level $\ell$    with process group index $\vecj$\\
    \midrule    $p^\ell_\vecj$, $p^\ell$ &    The process group at level $\ell$ with/without index $\vecj$\\
    \bottomrule
  \end{tabular}
  \caption{Commonly used notations.}
  \label{tab:Notations}
\end{table}

\subsection{Sparse Elimination}
\label{sub:Sparse Elimination}

Suppose that $A$ is a symmetric matrix. The row/column indices of $A$ are
partitioned into three sets $I \bigcup F \bigcup R$ where $I$ refers to the
interior point set, $F$ refers to the surrounding face point set, and $R$
refers to the rest point set. We further assume that there is no interaction
between the indices in $I$ and the ones in $R$. As a result, one can write
$A$ in the following form
\begin{equation}\label{eq:sparse elimination original A}
  A =
  \begin{bmatrix}
    A_{II} & A_{FI}^T & \\
    A_{FI} & A_{FF} & A_{RF}^T \\
    & A_{RF} & A_{RR} \\
  \end{bmatrix}.
\end{equation}

Let the $LDL^T$ decomposition of $A_{II}$ be $A_{II} = L_I D_I L_I^T$,
where $L_I$ is lower triangular matrix with unit diagonal. According
to the block Gaussian elimination of $A$ given by \eqref{eq:sparse
  elimination original A}, one defines the {\em sparse elimination} to
be
\begin{equation}\label{eq:SAS}
  S_I^T A S_I =
  \begin{bmatrix}
    D_I & &\\
    & B_{FF} & A_{RF}^T \\
    & A_{RF} & A_{RR}\\
  \end{bmatrix},
\end{equation}
where $B_{FF} = A_{FF}-A_{FI}A_{II}^{-1}A_{FI}^T$ is the associated
Schur complement and the explicit expressions for $S_I$ is
\begin{equation}
  S_I =
  \begin{bmatrix}
    L_{I}^{-T} & -A_{II}^{-1}A_{FI}^T & \\
    & I & \\
    & & I \\
  \end{bmatrix}.
\end{equation}

The sparse elimination removes the interaction between the interior
points $I$ and the corresponding surrounding face points $F$ and
leaves $A_{RF}$ and $A_{RR}$ untouched. We call the entire point set,
$I\bigcup F\bigcup R$, the {\em active point set}.  Then, after the
sparse elimination, the interior points are decoupled from other
points, which is conceptually equivalent to eliminate the interior
points from the active point set. After this, the new active point set
can be regarded as $F\bigcup R$.

\begin{figure}[htp]
    \centering
    \includegraphics[height=1.5in]{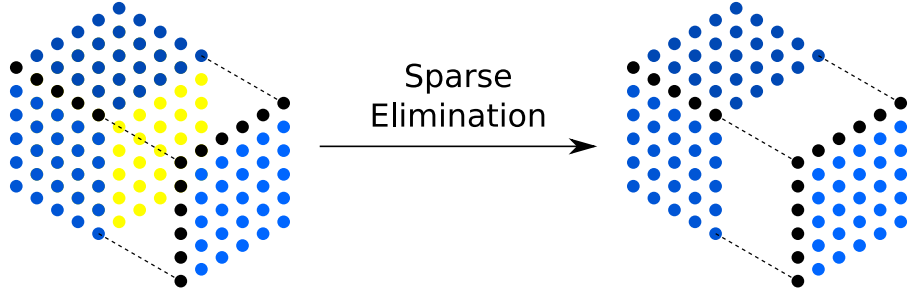}
    \caption{Sparse elimination:
    the interior points are eliminated after the sparse elimination;
    the rest points are not all plotted in the figure.}
    \label{fig:sparseelim}
\end{figure}

Figure~\ref{fig:sparseelim} illustrates the impact of the sparse
elimination. The dots in the figure represent the active points.
Before the sparse elimination (left), edge points, face points and
interior points are active while after the sparse elimination (right)
the interior points are eliminated from the active point set.

\subsection{Skeletonization}
\label{sub:Skeletonization}

Skeletonization is a tool for eliminating redundant point set from a
symmetric matrix that has low-rank off-diagonal blocks. The key step
in skeletonization uses the interpolative decomposition
\cite{Cheng2005,Martinsson2015} of low-rank matrices.

Let $A$ be a symmetric matrix of the form,
\begin{equation}
    A =
    \begin{bmatrix}
        A_{FF} & A_{RF}^T\\
        A_{RF} & A_{RR}
    \end{bmatrix},
\end{equation}
where $A_{RF}$ is a numerically low-rank matrix.  The interpolative
decomposition of $A_{RF}$ is (up to a permutation)
\begin{equation}
  A_{RF} =
  \begin{bmatrix}
    A_{R\cF} & A_{R\hF}\\
  \end{bmatrix}
  \approx
  \begin{bmatrix}
    A_{R\hF}T_F & A_{R\hF}\\
  \end{bmatrix},
\end{equation}
where $T_F$ is the interpolation matrix, $\hF$ is the skeleton point
set, $\cF$ is the redundant point set, and $F = \hF \bigcup \cF$.
Applying this approximation to $A$ results
\begin{equation}
  A \approx
  \left[
  \begin{array}{cc|c}
    A_{\cF\cF} & A_{\hF\cF}^T & T_F^TA_{R\hF}^T\\
    A_{\hF\cF} & A_{\hF\hF} & A_{R\hF}^T \bigstrut[b] \\
    \hline \bigstrut[t]
    A_{R\hF}T_F & A_{R\hF} & A_{RR}
  \end{array}\right],
\end{equation}
and be symmetrically factorized as
\begin{equation}\label{eq:SQAQS}
    S_{\cF}^T Q_{F}^{T} A Q_F S_{\cF} \approx
    S_{\cF}^T
    \left[
    \begin{array}{cc|c}
        B_{\cF\cF} & B_{\hF\cF}^T & \\
        B_{\hF\cF} & A_{\hF\hF} & A_{R\hF}^T \bigstrut[b]\\
        \hline \bigstrut[t]
        & A_{R\hF} & A_{RR}
    \end{array} \right]
    S_{\cF}
    =
    \left[ \begin{array}{cc|c}
         D_{\cF} & & \\
         & B_{\hF\hF} & A_{R\hF}^T \bigstrut[b]\\
         \hline \bigstrut[t]
         & A_{R\hF} & A_{RR}
    \end{array} \right],
\end{equation}
where
\begin{eqnarray}
    B_{\cF\cF} &=& A_{\cF\cF} - T_F^T A_{\hF\cF} - A_{\hF\cF}^T T_F
                + T_F^T A_{\hF\hF} T_F,\\
    B_{\hF\cF} &=& A_{\hF\cF} - A_{\hF\hF}T_F,\\
    B_{\hF\hF} &=& A_{\hF\hF} - B_{\hF\cF} B_{\cF\cF}^{-1} B_{\hF\cF}^T.
\end{eqnarray}
The factor $Q_F$ is generated by the block Gaussian elimination, which
is defined to be 
\begin{equation}
    Q_F =
    \left[ \begin{array}{cc|c}
        I & & \\
        -T_F& I &  \bigstrut[b]\\
        \hline \bigstrut[t]
        & & I \\
    \end{array} \right].
\end{equation}
Meanwhile, the factor $S_{\cF}$ is introduced in the sparse
elimination:
\begin{equation}
    S_{\cF} =
    \left[ \begin{array}{cc|c}
      L_{\cF}^{-T} & -B_{\cF\cF}^{-1}B_{\hF\cF}^T & \\
      & I &  \bigstrut[b]\\
      \hline \bigstrut[t]
      & & I \\
      \end{array} \right]
\end{equation}
where $L_{\cF}$ and $D_{\cF}$ come from the $LDL^T$ factorization of
$B_{\cF\cF}$, i.e., $B_{\cF\cF} = L_{\cF}D_{\cF}L_{\cF}^T$. Similar to
what happens in Section~\ref{sub:Sparse Elimination}, the
skeletonization eliminates the redundant point set $\cF$ from the
active point set.

\begin{figure}[htp]
    \centering
    \includegraphics[height=1.5in]{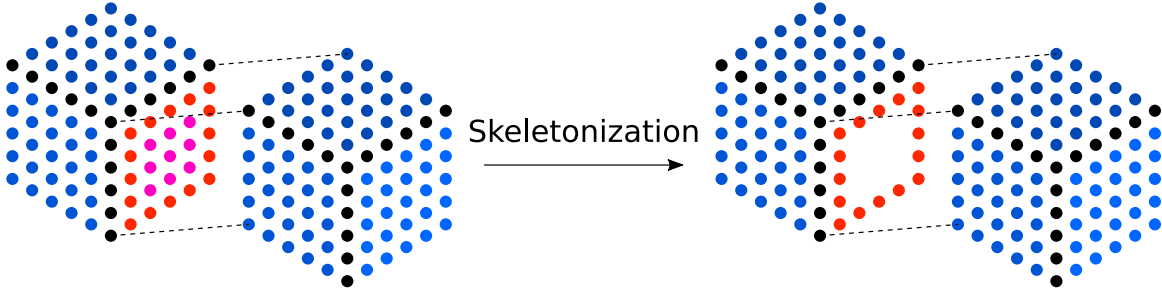}
    \caption{Skeletonization:
      the working face is colored by red and pink;
      red points are the skeleton points on the face
      whereas pink points are the redundant points on the face;
      skeletonization eliminates the redundant points
      from the active point set.}
    \label{fig:skel}
\end{figure}

The point elimination idea of the skeletonization is illustrated in
Figure~\ref{fig:skel}. Before the skeletonization (left), the edge
points, interior points, skeleton face points (red) and redundant face
points (pink) are all active, while after the skeletonization (right)
the redundant face points are eliminated from the active point set.

\subsection{Sequential HIF}
\label{sub:Sequential HIF}

This section reviews the sequential hierarchical interpolative
factorization (HIF) for 3D elliptic problems \eqref{eq:original-prob}
with the periodic boundary condition. Without loss of generality, we
discretize \eqref{eq:original-prob} with the seven-point stencil on a
uniform grid, which is defined in Section~\ref{sub:Notations}.  The
discrete system is
\begin{equation}\label{eq:seven-point stencil}
  \begin{split}
    \frac{1}{h^2}\left(a_{\vecj-\frac{1}{2}\vece_1} + a_{\vecj+\frac{1}{2}\vece_1}
    + a_{\vecj-\frac{1}{2}\vece_2} + a_{\vecj+\frac{1}{2}\vece_2}
    + a_{\vecj-\frac{1}{2}\vece_3} + a_{\vecj+\frac{1}{2}\vece_3}\right) u_\vecj - & \\
    \frac{1}{h^2}\left(
    a_{\vecj-\frac{1}{2}\vece_1} u_{\vecj-\vece_1}
    + a_{\vecj+\frac{1}{2}\vece_1} u_{\vecj+\vece_1}
    + a_{\vecj-\frac{1}{2}\vece_2} u_{\vecj-\vece_2}
    + a_{\vecj+\frac{1}{2}\vece_2} u_{\vecj+\vece_2}
    + a_{\vecj-\frac{1}{2}\vece_3} u_{\vecj-\vece_3}
    + a_{\vecj+\frac{1}{2}\vece_3} u_{\vecj+\vece_3}
    \right) & \\
    + b_\vecj u_\vecj = f_\vecj &
    \end{split}
\end{equation}
at each grid point $\vecx_\vecj$ for $\vecj = (j_1,j_2,j_3)$ and
$0\leq j_1,j_2,j_3 < n$,
where $a_\vecj = a(\vecx_\vecj)$, $b_\vecj = b(\vecx_\vecj)$,
$f_\vecj = f(\vecx_\vecj)$,
and $u_\vecj$ approximates the unknown function $u(x)$ at $\vecx_\vecj$.
The corresponding linear system is
\begin{equation}\label{eq:discretized linear system}
    Au = f
\end{equation}
where $A$ is a sparse SPD matrix if $b> 0$.

We first introduce the notion of active and inactive DOFs.
\begin{itemize}
\item A set $\Sigma$ of DOFs of $A$ are called {\bf active} if
  $A_{\Sigma\Sigma}$ is not a diagonal matrix or
  $A_{\bar{\Sigma}\Sigma}$ is a non-zero matrix;
\item A set $\Sigma$ of DOFs of $A$ are called {\bf inactive} if
  $A_{\Sigma\Sigma}$ is a diagonal matrix and $A_{\bar{\Sigma}\Sigma}$
  is a zero matrix.
\end{itemize}
Here $\bar{\Sigma}$ refers to the complement of the set $\Sigma$.
Sparse elimination and skeletonization provide concrete examples of
active and inactive DOFs. For example, sparse elimination turns the
indices $I$ from active DOFs of $A$ to inactive DOFs of $\wtA=S_I^T A
S_I$ in \eqref{eq:SAS}. Skeletonization turns the indices $\cF$ from
active DOFs of $A$ to inactive DOFs of $\wtA = S_{\cF}^T Q_F^T A
Q_FS_{\cF}$ in \eqref{eq:SQAQS}.

With these notations, the sequential HIF in \cite{Ho2015} is
summarized as follows. A more illustrative representation of the
sequential HIF is given on the left column of Figure~\ref{fig:DHIF}.

\begin{itemize}
\item {\bf Preliminary.} Let $A^0=A$ be the sparse symmetric matrix
  in \eqref{eq:discretized linear system}, $\Sigma^0$ be the initial
  active DOFs of $A$, which includes all indices.
  
\item {\bf Level $\ell$ for $\ell=0,\ldots,L-1$.}
  \begin{itemize}
  \item {\bf Preliminary.}  Let $A^\ell$ denote the matrix before any
    elimination at level $\ell$. $\Sigma^\ell$ is the corresponding
    active DOFs. Let us recall the notations in
    Section~\ref{sub:Notations}.  $C^\ell_\vecj$ denotes the active
    DOFs in the cell at level $\ell$ indexed with $\vecj$.
    $\calF^\ell_\vecj$ and $I^\ell_\vecj$ denote the surrounding faces
    and interior active DOFs in the corresponding cell, respectively.
  \item {\bf Sparse Elimination.}  We first focus on a single cell at
    level $\ell$ indexed with $\vecj$, i.e., $C^\ell_\vecj$.  To
    simplify the notation, we drop the superscript and subscript for
    now and introduce $C=C^\ell_\vecj$, $I=I^\ell_\vecj$,
    $F=\calF^\ell_\vecj$, and $R=R^\ell_\vecj$.  Based on the
    discretization and previous level eliminations, the interior
    active DOFs interact only with itself and its surrounding faces.
    The interactions of the interior active DOFs and the rest DOFs are
    empty and the corresponding matrix is zero, $A^\ell(R,I) = 0$.
    Hence, by applying sparse elimination, we have,
    \begin{equation}
      S_I^T A^\ell S_I =
      \begin{bmatrix}
        D_I & &\\
        & B^\ell_{F F} & \wideT{A^\ell_{R F}} \\
        & A^\ell_{R F} & A^\ell_{R R}\\
      \end{bmatrix},
    \end{equation}
    where the explicit definitions of $B^\ell_{F F}$ and $S_{I}$ are
    given in the discussion of sparse elimination. This factorization
    eliminates $I$ from the active DOFs of $A^\ell$.
    
    Looping over all cells $C^\ell_\vecj$ at level $\ell$, we obtain
    \begin{eqnarray}
      \wtA^\ell &=& \left(\prod_{I\in\calI^\ell}S_{I}\right)^T
      A^\ell \left(\prod_{I\in\calI^\ell}S_{I}\right),\\
      \wts^\ell &=& \Sigma^\ell \setminus \bigcup_{I\in\calI^\ell} I.
    \end{eqnarray}
    Now all the active interior DOFs at level $\ell$ are eliminated
    from $\Sigma^\ell$.
  \item {\bf Skeletonization.}  Each face at level $\ell$ not only
    interacts within its own cell but also interacts with faces of
    neighbor cells. Since the interaction between any two different
    faces is low-rank, this leads us to apply skeletonization. The
    skeletonization for any face $F \in \calF^\ell$ gives,
    \begin{equation}\label{eq:sequential HIF skeletonization}
      S_{\cF}^{T} Q_{F}^{T} \wtA^\ell Q_{F} S_{\cF} =
      \left[
        \begin{array}{cc|c}
          \wtD_{\cF} & & \\
          & \wtB^\ell_{\hF\hF}
          & \wideT{\wtA^\ell_{R\hF}} \bigstrut[b]\\
          \hline \bigstrut[t]
          & \wtA^\ell_{R\hF}
          & \wtA^\ell_{R R}
        \end{array}\right],
    \end{equation}
    where $\hF$ is the skeleton DOFs of $F$, $\cF$ is the redundant
    DOFs of $F$, and $R$ refers to the rest DOFs.  Due to the
    elimination from previous levels, $|F|$ scales as $\O(m2^\ell)$
    and $\wtA^\ell_{R F}$ contains a non-zero submatrix of size
    $\O(m2^\ell)\times \O(m2^\ell)$.  Therefore, the interpolative
    decomposition can be formed efficiently.  Readers are referred to
    Section~\ref{sub:Skeletonization} for the explicit forms of each
    matrix in \eqref{eq:sequential HIF skeletonization}.
    
    Looping over all faces at level $\ell$, we obtain
    \begin{equation}
      \begin{split}
        A^{\ell+1} \approx&
        \left(\prod_{F\in\calF^\ell}S_{\cF} Q_{F}\right)^T
        \wtA^\ell
        \left(\prod_{F\in\calF^\ell}S_{\cF} Q_{F}\right)\\
        = & \left(\prod_{F\in\calF^\ell}S_{\cF} Q_{F}\right)^T
        \left(\prod_{I\in\calI^\ell}S_{I}\right)^T
        A^\ell \left(\prod_{I\in\calI^\ell}S_{I}\right)
        \left(\prod_{F\in\calF^\ell}S_{\cF} Q_{F}\right)\\
        = & \wideT{W^\ell} A^\ell W^\ell,\\
      \end{split}
    \end{equation}
    where $W^\ell = \left(\prod_{I\in\calI^\ell}S_{I}\right)
    \left(\prod_{F\in\calF^\ell}S_{\cF} Q_{F}\right)$. The active DOFs
    for the next level is now defined as,
    \begin{equation}
      \Sigma^{\ell+1} = \wts^\ell \setminus \bigcup_{F\in\calF^\ell}
      \cF = \Sigma^{\ell} \setminus \left(
      \left(\bigcup_{I\in\calI^\ell} I\right) \bigcup
      \left(\bigcup_{F\in\calF^\ell} \cF \right) \right).
    \end{equation}
  \end{itemize}

\item {\bf Level $L$.} Finally, $A^L$ and $\Sigma^L$ are the matrix
  and active DOFs at level $L$. Up to a permutation, $A^L$ can be
  factorized as
  \begin{equation}
    A^L =
    \begin{bmatrix}
      A^L_{\Sigma^L\Sigma^L} & \\
      & D_R
    \end{bmatrix}
    =
    \begin{bmatrix}
      L_{\Sigma^L} & \\
      & I
    \end{bmatrix}
    \begin{bmatrix}
      D_{\Sigma^L} & \\
      & D_R
    \end{bmatrix}
    \begin{bmatrix}
      L_{\Sigma^L}^T & \\
      & I
    \end{bmatrix}
    := \wideInvT{W^L}D\wideInv{W^L}.
  \end{equation}
  
  Combining all these factorization results
  \begin{equation}
    A \approx \wideInvT{W^0}\cdots\wideInvT{W^{L-1}}\wideInvT{W^L}
    D\wideInv{W^L}\wideInv{W^{L-1}}\cdots\wideInv{W^0} \defeq F
  \end{equation}
  and
  \begin{equation}
    A^{-1} \approx W^0\cdots W^{L-1} W^L D^{-1}
    \wideT{W^L}\wideT{W^{L-1}}\cdots\wideT{W^0} = F^{-1}.
  \end{equation}
  $A^{-1}$ is factorized into a multiplicative sequence of
  matrices $W^\ell$ and each $W^\ell$ corresponding to level $\ell$
  is again a multiplicative sequence of sparse matrices, $S_I$, $S_{\cF}$
  and $Q_F$. Due to the fact that any $S_I$, $S_{\cF}$ or $Q_F$ contains
  a small non-trivial (i.e., neither identity nor empty) matrix of size
  $\O(\frac{N^{1/3}}{2^{L-\ell}})\times \O(\frac{N^{1/3}}{2^{L-\ell}})$,
  the overall complexity for strong and applying $W^\ell$ is
  $\O(N/2^\ell)$. Hence the application of the inverse of $A$ is of
  $\O(N)$ computation and memory complexity.
\end{itemize}

\section{Distributed-memory HIF}
\label{sec:Distributed-memory HIF}

This section describes the algorithm for the distributed-memory HIF.

\subsection{Process Tree}
\label{sub:Process Tree}

\begin{figure}[htp]
    \centering
    \includegraphics[height=1.5in]{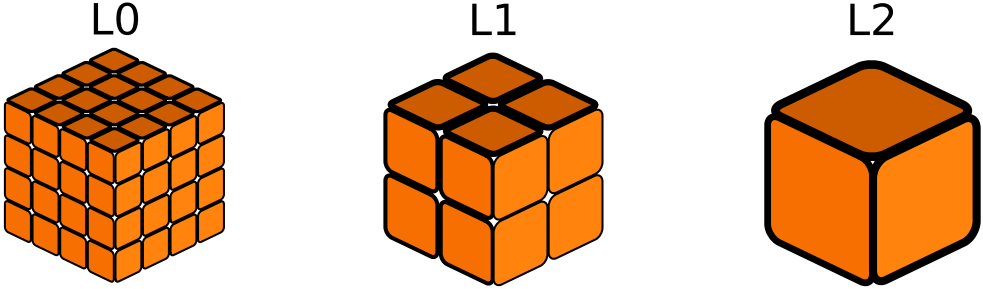}
    \caption{Process tree:
    64 processes are organized in the process tree.}
    \label{fig:processtree}
\end{figure}

For simplicity, assume that there are $8^L$ processes. We introduce a
{\em process tree} that has $L+1$ levels and resembles the
hierarchical structure of the computation domain.  Each node of this
process tree is called a {\em process group}. First at the leaf level,
there are $8^L$ leaf process groups denoted as
$\{p^0_\vecj\}_\vecj$. Here $\vecj=(j_1,j_2,j_3)$, $0\leq
j_1,j_2,j_3<2^L$ and the superscript $0$ refers to the leaf level
(level 0). Each group at this level only contains a single process.
Each node at level 1 of the process tree is constructed by merging 8
leaf processes. More precisely, we denote the process group at level 1
as $p^1_\vecj$ for $\vecj=(j_1,j_2,j_3)$, $0\leq j_1,j_2,j_3 <
2^{L-1}$, and $p^1_\vecj = \bigcup_{\floor{\vecj_c/2} = \vecj}
p^0_{\vecj_c}$.  Similarly, we recursively define the node at level
$\ell$ as $p^\ell_\vecj = \bigcup_{\floor{\vecj_c/2} = \vecj}
p^{\ell-1}_{\vecj_c}$.  Finally, the process group $p^L_\vecZero$ at
the root includes all processes.  Figure~\ref{fig:processtree}
illustrates the process tree.  Each cube in the process tree is a
process group.

\subsection{Distributed-memory method}
\label{sub:Distributed-memory method}

Same as in Section~\ref{sub:Sequential HIF}, we define the $n \times
n\times n$ grid on $\Omega=[0,1)^3$ for $n=m2^L$, where $m=\O(1)$ is a
  small positive integer and $L=\O(\log N)$ is the level number of the
  root level.  Discretizing \eqref{eq:original-prob} with seven-point
  stencil on the grid provides the linear system $Au=f$, where $A$ is
  a sparse $N\times N$ SPD matrix, $u\in \bbR^N$ is the unknown
  function at grid points, and $f\in \bbR^N$ is the given function at
  grid points.

Given the process tree (Section~\ref{sub:Process Tree}) with $8^L$
processes and the sequential HIF structure
(Section~\ref{sub:Sequential HIF}), the construction of the
distributed-memory hierarchical interpolative factorization (DHIF)
consists of the following steps.

\begin{itemize}
\item {\bf Preliminary.}  Construct the process tree with $8^L$
  processes. Each process group $p^0_\vecj$ owns the data
  corresponding to cell $C^0_\vecj$ and the set of active DOFs in
  $C^0_\vecj$ are denoted as $\Sigma^0_\vecj$, for
  $\vecj=(j_1,j_2,j_3)$ and $0\leq j_1,j_2,j_3 < 2^L$. Set $A^0=A$
  and let the process group $p^0_\vecj$ own $A^0(:,\Sigma^0_\vecj)$,
  which is a sparse tall-skinny matrix with $\O(N/P)$ non-zero
  entries.
  
\item {\bf Level $\ell$ for $\ell=0,\ldots,L-1$.}
  \begin{itemize}
  \item {\bf Preliminary.}  Let $A^\ell$ denote the matrix before any
    elimination at level $\ell$.  $\Sigma^\ell_\vecj$ denotes the
    active DOFs owned by the process group $p^\ell_\vecj$ for
    $\vecj=(j_1,j_2,j_3)$, $0\leq j_1,j_2,j_3 < 2^{L-\ell}$, and the
    non-zero submatrix of $A^\ell(:,\Sigma^\ell_\vecj)$ is distributed
    among the process group $p^\ell_\vecj$ using the two-dimensional
    block-cyclic distribution. 
  \item {\bf Sparse Elimination.} The process group $p^\ell_\vecj$
    owns $A^\ell(:,\Sigma^\ell_\vecj)$, which is sufficient for
    performing sparse elimination for $I^\ell_\vecj$. To simplify the
    notation, we define $I=I^\ell_\vecj$ as the active interior DOFs
    of cell $C^\ell_\vecj$, $F=\calF^\ell_\vecj$ as the surrounding
    faces, and $R=R^\ell_\vecj$ as the rest active DOFs. Sparse
    elimination at level $\ell$ within the process group
    $p^\ell_\vecj$ performs essentially
    \begin{equation}
      S_I^{T} A^\ell S_I =
      \begin{bmatrix}
        D_I & &\\
        & B^\ell_{FF} & \wideT{A^\ell_{RF}}\\
        & A^\ell_{RF} & A^\ell_{RR}
      \end{bmatrix},
    \end{equation}
    where $B^\ell_{FF} = A^\ell_{FF}-A^\ell_{FI}\wideInv{A^\ell_{II}}
    \wideT{A^\ell_{FI}}$,
    \begin{equation}
      S_I =
      \begin{bmatrix}
          \wideInvT{L^\ell_I} &
          -\wideInv{A^\ell_{II}}\wideT{A^\ell_{FI}} &\\
          & I &\\
          & & I
        \end{bmatrix}
    \end{equation}
    with $L^\ell_I D_I \wideT{L^\ell_I} = A^\ell_{II}$.  Since
    $A^\ell(:,\Sigma^\ell_\vecj)$ is owned locally by $p^\ell_\vecj$,
    both $A^\ell_{FI}$ and $A^\ell_{II}$ are local matrices.  All
    non-trivial (i.e., neither identity nor empty) submatrices in
    $S_I$ are formed locally and stored locally for application.  On
    the other hand, updating on $A^\ell_{FF}$ requires some
    communication in the next step.
  \item {\bf Communication after sparse elimination.}  After all
    sparse eliminations are performed, some communication is required
    to update $A^\ell_{FF}$ for each cell $C^\ell_\vecj$.  For the
    problem \eqref{eq:original-prob} with the periodic boundary
    conditions, each face at level $\ell$ is the surrounding face of
    exactly two cells. The owning process groups of these two cells
    need to communicate with each other to apply the additive updates,
    a submatrix of
    $-A^\ell_{FI}\wideInv{A^\ell_{II}}\wideT{A^\ell_{FI}}$.  Once all
    communications are finished, the parallel sparse elimination does
    the rest of the computation, which can be conceptually denoted as,
    \begin{equation}
      \begin{split}
        & \wtA^\ell = \wideT{\prod_{I\in\calI^\ell}S_{I}}
        A^\ell \left(\prod_{I\in\calI^\ell}S_{I}\right),\\
        & \wts^\ell_\vecj = \Sigma^\ell_\vecj \setminus
        \bigcup_{I\in\calI^\ell}I,
      \end{split}
    \end{equation}
    for $\vecj=(j_1,j_2,j_3), 0\leq j_1,j_2,j_3<2^{L-\ell}$.
    
  \item {\bf Skeletonization.}  For each face $F$ owned by
    $p^\ell_\vecj$, the corresponding matrices $\wtA^\ell(:,F)$ is
    stored locally. Similar to the parallel sparse elimination part,
    most operations are local at the process group $p^\ell_\vecj$ and
    can be carried out using the dense parallel linear algebra
    efficiently. By forming a parallel interpolative decomposition
    (ID) for $\wtA^\ell_{RF} =
    \begin{bmatrix} \wtA^\ell_{R\hF} T^\ell_F &
    \wtA^\ell_{R\hF} \end{bmatrix}$, the parallel skeletonization can
    be, conceptually, written as,
    \begin{equation}
      S_{\cF} Q_F
      \wtA^\ell \wideT{Q_F} \wideT{S_{\cF}} \approx
      \left[\begin{array}{cc|c}
          D_{\cF} & & \\
          & \wtB^\ell_{\hF\hF} & \wtA^\ell_{R\hF} \bigstrut[b]\\
          \hline \bigstrut[t]
          & \wtA^\ell_{R\hF} & \wtA^\ell_{RR}\\
        \end{array} \right],
    \end{equation}
    where the definitions of $Q_F$ and $S_{\cF}$ are given in the
    discussion of skeletonization. Since $\wtA^\ell_{\cF\cF}$,
    $\wtA^\ell_{\hF\cF}$, $\wtA^\ell_{\hF\hF}$ and $T^\ell_{F}$ are
    all owned by $p^\ell_\vecj$, it requires only local operations to
    form
    \begin{equation}
      \begin{split}
        \wtB^\ell_{\cF\cF} = &
        \wtA^\ell_{\cF\cF} - \wideT{T^\ell_F}\wtA^\ell_{\hF\cF}
        - \wideT{\wtA^\ell_{\hF\cF}}T^\ell_F
        + \wideT{T^\ell_F}\wtA^\ell_{\hF\hF}T^\ell_F,\\
        \wtB^\ell_{\hF\cF} = & \wtA^\ell_{\hF\cF}
        - \wtA^\ell_{\hF\hF}T^\ell_F,\\
        \wtB^\ell_{\hF\hF} = &
        \wtA^\ell_{\hF\hF}
        - \wtB^\ell_{\hF\cF}\wideInv{\wtB^\ell_{\cF\cF}}
        \wideT{\wtB^\ell_{\hF\cF}}.\\
      \end{split}
    \end{equation}
    Similarly, $L_{\cF}$, which is the $LDL^T$ factor of
    $\wtB_{\cF\cF}$, is also formed within the process group
    $p^\ell_\vecj$. Moreover, since non-trivial blocks in $Q_F$ and
    $S_{\cF}$ are both local, this implies that the applications of
    $Q_F$ and $S_{\cF}$ are local operations. As a result, the
    parallel skeletonization factorizes $A^\ell$ conceptually as,
    \begin{equation}
      \begin{split}
        A^{\ell+1} \approx &
        \wideT{\prod_{F\in\calF^\ell}S_{\cF}Q_{F}}
        \wtA^\ell
        \left(\prod_{F\in\calF^\ell}S_{\cF}Q_{F}\right)\\
        = &
        \wideT{\prod_{F\in\calF^\ell}S_{\cF}Q_{F}}
        \wideT{\prod_{I\in\calI^\ell}S_{I}}
        A^\ell
        \left(\prod_{I\in\calI^\ell}S_{I}\right)
        \left(\prod_{F\in\calF^\ell}S_{\cF}Q_{F}\right)\\
      \end{split}
    \end{equation}
    and we can define
    \begin{equation}
      \begin{split}
        W^{\ell} = & \left(\prod_{I\in\calI^\ell}S_{I}\right)
        \left(\prod_{F\in\calF^\ell}S_{\cF}Q_{F}\right),\\
        \Sigma^{\ell+1/2}_\vecj = &
        \wts^\ell_\vecj \setminus \bigcup_{F\in\calF^\ell}\cF\\
        = & \Sigma^\ell_\vecj \setminus
        \left(\left(\bigcup_{F\in\calF^\ell}\cF\right)\bigcup
        \left(\bigcup_{I\in\calI^\ell}I\right)\right).
      \end{split}
    \end{equation}
    We would like to emphasize that the factors $W^{\ell}$ are
    evenly distributed among the process groups at level $\ell$ and
    that all non-trivial blocks are stored locally.
    
  \item {\bf Merging and Redistribution.}  Towards the end of the
    factorization at level $\ell$, we need to merge the process groups
    and redistribute the data associated with the active DOFs in order
    to prepare for the work at level $\ell+1$.  For each process group
    at level $\ell+1$, $p^{\ell+1}_\vecj$, for $\vecj=(j_1,j_2,j_3)$,
    $0\leq j_1,j_2,j_3 < 2^{L-\ell-1}$, we first form its active DOF
    set $\Sigma^{\ell+1}_\vecj$ by merging
    $\Sigma^{\ell+1/2}_{\vecj_c}$ from all its children
    $p^\ell_{\vecj_c}$, where $\floor{\vecj_c/2} = \vecj$. In
    addition, $A^{\ell+1}(:,s^{\ell+1}_\vecj)$ is separately owned by
    $\left\{p^\ell_{\vecj_c}\right\}_{\floor{\vecj_c/2}=\vecj}$. A
    redistribution among $p^{\ell+1}_\vecj$ is needed in order to
    reduce the communication cost for future parallel dense linear
    algebra. Although this redistribution requires a global
    communication among $p^{\ell+1}_\vecj$, the complexities for
    message and bandwidth are bounded by the cost for parallel dense
    linear algebra. Actually, as we shall see in the numerical
    results, its cost is far lower than that of the parallel dense
    linear algebra.
  \end{itemize}

\item {\bf Level $L$ Factorization.}  The parallel factorization at
  level $L$ is quite similar to the sequential one.  After
  factorizations from all previous levels,
  $A^L(\Sigma^L_\vecZero,\Sigma^L_\vecZero)$ is distributed among
  $p^L_\vecZero$.  A parallel $LDL^T$ factorization of
  $A^L_{\Sigma^L_\vecZero \Sigma^L_\vecZero} =
  A^L(\Sigma^L_\vecZero,\Sigma^L_\vecZero)$ among the processes in
  $p^L_\vecZero$ results
  \begin{equation}
    A^L =
    \begin{bmatrix}
      A^L_{\Sigma^L_\vecZero \Sigma^L_\vecZero} & \\
      & D_R
    \end{bmatrix}
    =
    \begin{bmatrix}
      L^L_{\Sigma^L_\vecZero} & \\
      & I
    \end{bmatrix}
    \begin{bmatrix}
      D^L_{\Sigma^L_\vecZero} & \\
      & D_R
    \end{bmatrix}
    \begin{bmatrix}
      \wideT{L^L_{\Sigma^L_\vecZero}} & \\
      & I
    \end{bmatrix}
    = \wideInvT{W^L}D\wideInv{W^L}.
  \end{equation}
  
  Consequently, we forms the DHIF for $A$ and $A^{-1}$ as
  \begin{equation}
    A \approx \wideInvT{W^0} \cdots \wideInvT{W^{L-1}} \wideInvT{W^L} D
    \wideInv{W^L} \wideInv{W^{L-1}} \cdots \wideInv{W^0} \defeq F
  \end{equation}
  and
  \begin{equation}
    A^{-1} \approx W^0 \cdots W^{L-1} W^L
    D^{-1} \wideT{W^L} \wideT{W^{L-1}} \cdots \wideT{W^0} = F^{-1}.
  \end{equation}
  
  The factors, $W^\ell$ are evenly distributed among all processes and
  the application of $F^{-1}$ is basically a sequence of parallel
  dense matrix-vector multiplications.
\end{itemize}

\begin{figure}[ht!]
  \centering
  \includegraphics[height=0.93\textheight]{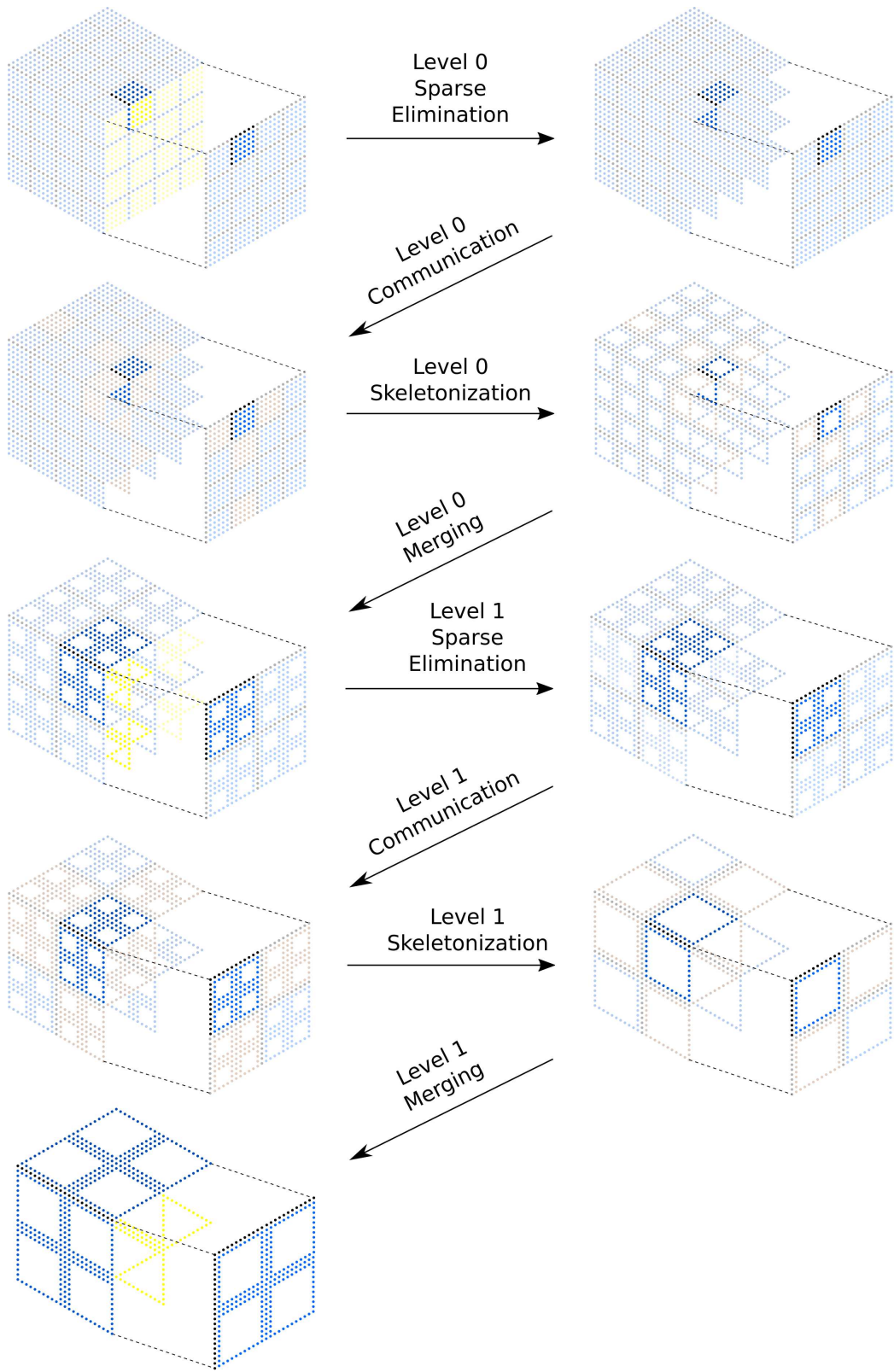}
  \caption{Distributed-memory hierarchical interpolative factorization.}
  \label{fig:DHIF}
\end{figure}

In Figure~\ref{fig:DHIF}, we illustrate an example of DHIF for problem
of size $24 \times 24 \times 24$ with $m=6$ and $L=2$. The computation
is distributed on a process tree with $4^3=64$ processes.
Particularly, Figure~\ref{fig:DHIF} highlights the DOFs owned by
process groups involving $p^0_{(0,1,0)}$, i.e., $p^0_{(0,1,0)}$,
$p^1_{(0,0,0)}$, and $p^2_{(0,0,0)}$.  Yellow points denote interior
active DOFs, blue and brown points denote face active DOFs, and black
points denote edge active DOFs.  Meanwhile, we also have unfaded and
faded groups of points.  Unfaded points are owned by the process
groups involving $p^0_{(0,1,0)}$.  In other words, $p^0_{(0,1,0)}$ is
the owner for part of the unfaded points.  The faded points are owned
by other process groups.  In the second row and the forth row, we also
see faded brown points, which indicates the required communication to
process $p^0_{(0,1,0)}$. Here Figure~\ref{fig:DHIF} works through two
levels of the elimination processes of the DHIF step by step.

\subsection{Complexity Analysis}
\label{sub:Complexity Analysis}

\subsubsection{Memory Complexity}
\label{ssub:Memory Complexity}

There are two places in the distributed algorithm that require heavy
memory usage.  The first one is to store the original matrix $A$ and
its updated version $A^\ell$ for each level $\ell$. As we mentioned
above in the parallel algorithm, $A^\ell$ contains at most $\O(N)$
non-zeros and they are evenly distributed on $P$ processes as
follows. At level
$\ell$, there are $8^{L-\ell}$ cells, and empirically each of which
contains $\O\left(\frac{N^{1/3}}{2^{L-\ell}}\right)$ active DOFs.
Meanwhile, each cell is evenly owned by a process group with $8^\ell$
processes.  Hence,
$\O\left(\left(\frac{N^{1/3}}{2^{L-\ell}}\right)^2\right)$ non-zero
entries of $A^\ell(:,s^\ell_\vecj)$ is evenly distributed on process
group $p^\ell_\vecj$ with $8^\ell$ processes.  Overall, there are
$\O\left(8^{L-\ell}\cdot \frac{N^{2/3}}{4^{L-\ell}}\right) =\O(N\cdot
2^{-\ell})$ non-zero entries in $A^\ell$ evenly distributed on
$8^{L-\ell}\cdot 8^\ell = P$ processes, and each process owns
$\O\left(\frac{N}{P}\cdot 2^{-\ell}\right)$ data for $A^\ell$.
Moreover, the factorization at level $\ell$ does not rely on the
matrix $A^{\ell'}$ for $\ell'<\ell-1$.  Therefore, the memory cost for
storing $A^\ell$s is $\O(\frac{N}{P})$ for each process.

The second place is to store the factors $W^\ell$.  It is not
difficult to see that the memory cost for each $W^\ell$ is the same as
$A^{\ell}$. Only non-trivial blocks in $S_I$, $Q_F$, and
$S_{\cF}$ require storage. Since each of these non-trivial blocks
is of size
$\O\left(\frac{N^{1/3}}{2^{L-\ell}}\right)\times\O\left(\frac{N^{1/3}}{2^{L-\ell}}\right)$
and evenly distributed on $8^\ell$ processes, the overall memory
requirement for each $W^\ell$ on a process is
$\O\left(\frac{N}{P}\cdot 2^{-\ell}\right)$. Therefore,
$\O\left(\frac{N}{P}\right)$ memory is required on each process to
store all $W^\ell$s.

\subsubsection{Computation Complexity}
\label{ssub:Computational Complexity}

The majority of the computation work goes to the construction of
$S_I$, $Q_{F}$ and $S_{\cF}$.  As stated in the
previous section, at level $\ell$, each non-trivial dense matrix in
these factors is of size
$\O\left(\frac{N^{1/3}}{2^{L-\ell}}\right)\times
\O\left(\frac{N^{1/3}}{2^{L-\ell}}\right)$.  The construction adopts
the matrix-matrix multiplication, the interpolative decomposition
(pivoting QR), the $LDL^T$ factorization, and the triangular matrix
inversion. Each of these operation is of cubic computation
complexities and the corresponding parallel computation cost over
$8^\ell$ processes is $\O\left(\frac{N}{P}\right)$. Since there is
only a constant number of these operations per process at a single
level, the total computational complexity across all
$\O(\log N)$ levels is $\O\left(\frac{N\log N}{P}\right)$.

The application computational complexity is simply the complexity of
applying each non-zero entries in $W^\ell$s once, hence, the overall
computational complexity is the same as the memory complexity
$\O\left(\frac{N}{P}\right)$.

\subsubsection{Communication Complexity}
\label{ssub:Communication Complexity}

The communication complexity is composed of three parts: the
communication in the parallel dense linear algebra, the communication
after sparse elimination, and the merging and redistribution step
within DHIF. It is clear to see that the communication cost for the
second part is bounded by either of the rest. Hence, we will simply
derive the communication cost for the first and third parts. Here, we
adopt the simplified communication model, $T_{comm} = \alpha + \beta$,
where $\alpha$ is the latency, and $\beta$ is the inverse bandwidth.

At level $\ell$, the parallel dense linear algebra involves the
matrix-matrix multiplication, the ID, the $LDL^T$ factorization, and
the triangular matrix inversion for matrices of size
$\O\left(\frac{N^{1/3}}{2^{L-\ell}}\right)\times\O\left(\frac{N^{1/3}}{2^{L-\ell}}\right)$.
All these basic operations are carried out on a process group of size
$8^\ell$. Following the discussion in \cite{Ballard2011}, the
communication cost for these operations are bounded by
$\O\left(\ell^3\sqrt{8^\ell}\right)\alpha+\O\left(\frac{N^{2/3}}{4^{L-\ell}8^\ell}\ell\right)\beta$.
Summing over all levels, one can control the communication cost of the
parallel dense linear algebra part by
\begin{equation}
    \O\left(\sqrt{P}\log^3 P \right)\alpha
    + \O\left(\frac{N^{2/3}}{P^{2/3}}\right)\beta.
\end{equation}

On the other hand, the merging and redistribution step at level $\ell$
involves $8^{\ell+1}$ processes and redistributes matrices of size
$\O\left(\frac{N^{1/3}}{2^{L-\ell}}\cdot
8\right)\times\O\left(\frac{N^{1/3}}{2^{L-\ell}}\cdot 8\right)$.  The
current implementation adopts the MPI routine \text{MPI\_AllToAll} to
handle the redistribution on a 2D process mesh.  Further, we assume
the all-to-all communication sends and receives long messages. The
standard upper bound for the cost of this routine is
$\O\left(\sqrt{8^{\ell+1}}\right)\alpha
+\O\left(\frac{N^{2/3}}{4^{L-\ell}\sqrt{8^{\ell+1}}}\cdot
64\right)\beta$~\cite{Scott1991}. Therefore, the
over all cost is
\begin{equation}
  \O\left(\sqrt{P}\right)\alpha + \O\left(\frac{N^{2/3}}{\sqrt{P}}\right)\beta.
\end{equation}
The complexity of the latency part is not scalable.  However
empirically, the cost for this communication is relatively small in
the actual running time.

\section{Numerical Results} \label{sec:Numerical Results}

Here we present a couple of numerical examples to demonstrate the parallel
efficiency of the distributed-memory HIF. The algorithm is implemented in
C++11 and all inter-process communication is expressed via the Message
Passing Interface (MPI). The distributed-memory dense linear algebra
computation is done through the Elemental library~\cite{Poulson2013}. All
numerical examples are performed on Edison at the National Energy Research
Scientific Computing center (NERSC). The numbers of processes used are
always powers of two, ranging from 1 to 8192. The memory allocated for
each process is limited to 2GB.

All numerical results are measured in two ways: the strong scaling and
weak scaling.  The strong scaling measurement fixes the problem size,
and increases the number of processes.  For a fixed problem size, let
$T^S_{P}$ be the running time of $P$ processes.  The strong scaling
efficiency is defined as,
\begin{equation}
  E^S = \frac{T^S_1}{P\cdot T^S_P}.
\end{equation}
In the case that $T^S_1$ is not available,
e.g., the fixed problem can not fit into the single process memory,
we adopts the first available running time, $T^S_m$, associating
with the smallest number of processes, $m$,
as a reference.
And the modified strong scaling efficiency is,
\begin{equation}
  E^S = \frac{m\cdot T^S_m}{P\cdot T^S_P}.
\end{equation}

The weak scaling measurement fixes the ratio between the problem size
and the number of processes. For a fixed ratio, we define the weak
scaling efficiency as,
\begin{equation}
  E^W = \frac{T^W_m}{T^W_P},
\end{equation}
where $T^W_m$ is the first available running time with $m$ processes,
and $T^W_P$ is the running time of $P$ processes.  

\begin{table}[htp]
  \centering
  \begin{tabular}{cm{10cm}}
    \toprule
    Notation & Explanation \\
    \toprule
    $\epsilon$ & Relative precision of the ID\\
    $N$ & Total number of DOFs in the problem\\
    $e_s$ & Relative error for solving, $\norm{(I-F^{-1}A)x}/\norm{x}$,
    where $x$ is a Gaussian random vector\\
    $|\Sigma_L|$ & Number of remaining active DOFs at the root level\\
    $m_f$ & Maximum memory required to perform the factorization in GB
    across all processes\\
    $t_f$ & Time for constructing the factorization in seconds\\
    $E^S_f$ & Strong scaling efficiency for factorization time\\
    $t_s$ & Time for applying $F^{-1}$ to a vector in seconds\\
    $E^S_s$ & Strong scaling efficiency for application time\\
    $n_{iter}$ & Number of iterations to solve $Au=f$ with GMRES
    with $F^{-1}$ being a preconditioner
    to a tolerance of $10^{-12}$\\
    \bottomrule
  \end{tabular}
  \caption{Notations for the numerical results}
  \label{tab:notations-num}
\end{table}

The notations used in the following tables and figures are explained
in Table~\ref{tab:notations-num}. For simplicity, all examples are
defined over $\Omega = [0,1)^3$ with periodic boundary condition,
  discretized on a uniform grid, $n\times n\times n$, with $n$ being
  the number of points in each dimension and $N=n^3$.  The PDEs
  defined in \eqref{eq:original-prob} is discretized using
  the second-order central difference method with seven-point stencil,
  which is the same as \eqref{eq:seven-point stencil}. Octrees are
  adopted to partition the computation domain with the block size at
  leaf level bounded by 64.
  
{\bf \emph{Example 1.}}
We first consider the problem in \eqref{eq:original-prob}
with $a(x) \equiv 1$ and $b(x) \equiv 0.1$.
The relative precision of the ID is set to be $\epsilon = 10^{-3}$.

\begin{table}
    \centering
    \begin{tabular}{crcrccrcrc}
\toprule
 $N$ & $P$ & $e_s$ & $|s_L|$ & $m_f$ & $t_f$ & $E^S_f$ & $t_s$ & $E^S_s$ & $n_{iter}$ \\
\toprule
\multirow{6}{*}{$32^3$}
& 1 & 4.84e-04 & 3440 & 1.92e-01 & 4.85e+00 & 100\% & 1.36e-01 & 100\% & 6 \\
& 2 & 5.26e-04 & 3440 & 9.60e-02 & 2.60e+00 &  93\% & 6.65e-02 & 103\% & 6 \\
& 4 & 3.78e-04 & 3440 & 4.80e-02 & 1.45e+00 &  84\% & 3.47e-02 &  98\% & 6 \\
& 8 & 4.93e-04 & 3440 & 2.40e-02 & 8.38e-01 &  72\% & 1.99e-02 &  85\% & 6 \\
& 16 & 3.97e-04 & 3440 & 1.20e-02 & 5.83e-01 &  52\% & 1.31e-02 &  65\% & 6 \\
& 32 & 7.33e-04 & 3440 & 6.03e-03 & 4.35e-01 &  35\% & 1.47e-02 &  29\% & 6 \\
\midrule
\multirow{8}{*}{$64^3$}
& 2 & 5.92e-04 & 7760 & 9.07e-01 & 3.87e+01 & 100\% & 6.08e-01 & 100\% & 6 \\
& 4 & 5.98e-04 & 7760 & 4.54e-01 & 2.36e+01 &  82\% & 2.99e-01 & 102\% & 6 \\
& 8 & 5.59e-04 & 7760 & 2.27e-01 & 1.48e+01 &  65\% & 1.61e-01 &  94\% & 6 \\
& 16 & 6.30e-04 & 7760 & 1.13e-01 & 1.03e+01 &  47\% & 9.52e-02 &  80\% & 6 \\
& 32 & 5.89e-04 & 7760 & 5.68e-02 & 5.34e+00 &  45\% & 6.88e-02 &  55\% & 6 \\
& 64 & 5.45e-04 & 7760 & 2.84e-02 & 2.67e+00 &  45\% & 4.10e-02 &  46\% & 6 \\
& 128 & 5.43e-04 & 7760 & 1.42e-02 & 1.52e+00 &  40\% & 3.43e-02 &  28\% & 6 \\
& 256 & 6.29e-04 & 7760 & 7.14e-03 & 1.27e+00 &  24\% & 2.69e-02 &  18\% & 6 \\
\midrule
\multirow{7}{*}{$128^3$}
& 16 & 6.19e-04 & 16208 & 9.77e-01 & 1.43e+02 & 100\% & 8.24e-01 & 100\% & 6 \\
& 32 & 5.98e-04 & 16208 & 4.89e-01 & 7.40e+01 &  97\% & 4.37e-01 &  94\% & 6 \\
& 64 & 5.85e-04 & 16208 & 2.44e-01 & 3.87e+01 &  92\% & 2.26e-01 &  91\% & 6 \\
& 128 & 6.23e-04 & 16208 & 1.22e-01 & 2.11e+01 &  85\% & 1.40e-01 &  74\% & 6 \\
& 256 & 6.14e-04 & 16208 & 6.12e-02 & 1.00e+01 &  89\% & 9.76e-02 &  53\% & 6 \\
& 512 & 5.96e-04 & 16208 & 3.06e-02 & 5.80e+00 &  77\% & 1.98e-01 &  13\% & 6 \\
& 1024 & 5.86e-04 & 16208 & 1.54e-02 & 3.46e+00 &  65\% & 6.13e-02 &  21\% & 6 \\
\midrule
\multirow{7}{*}{$256^3$}
& 128 & 6.18e-04 & 33104 & 1.01e+00 & 2.24e+02 & 100\% & 9.18e-01 & 100\% & 6 \\
& 256 & 6.11e-04 & 33104 & 5.07e-01 & 1.19e+02 &  94\% & 4.88e-01 &  94\% & 6 \\
& 512 & 6.06e-04 & 33104 & 2.53e-01 & 6.33e+01 &  88\% & 2.85e-01 &  81\% & 6 \\
& 1024 & 6.25e-04 & 33104 & 1.27e-01 & 3.19e+01 &  88\% & 1.86e-01 &  62\% & 6 \\
& 2048 & 6.18e-04 & 33104 & 6.35e-02 & 2.44e+01 &  57\% & 1.58e-01 &  36\% & 6 \\
& 4096 & 6.16e-04 & 33104 & 3.18e-02 & 1.27e+01 &  55\% & 1.73e-01 &  17\% & 6 \\
& 8192 & 6.14e-04 & 33104 & 1.60e-02 & 1.16e+01 &  30\% & 4.14e-01 &   3\% & 6 \\
\midrule
\multirow{4}{*}{$512^3$}
& 1024 & 6.16e-04 & 66896 & 1.03e+00 & 3.32e+02 & 100\% & 1.08e+00 & 100\% & 6 \\
& 2048 & 6.15e-04 & 66896 & 5.16e-01 & 1.84e+02 &  90\% & 6.53e-01 &  82\% & 6 \\
& 4096 & 6.14e-04 & 66896 & 2.58e-01 & 9.55e+01 &  87\% & 4.90e-01 &  55\% & 6 \\
& 8192 & 6.13e-04 & 66896 & 1.29e-01 & 5.58e+01 &  74\% & 4.58e-01 &  29\% & 6 \\
\midrule
\multirow{1}{*}{$1024^3$}
& 8192 & 6.15e-04 & 134480 & 1.04e+00 & 4.67e+02 & 100\% & 1.48e+00 & 100\% & 6 \\
\bottomrule
    \end{tabular}
    \caption{Example 1. Numerical results.}
    \label{tab:num-ex1}
\end{table}

\begin{figure*}[htp]
    \centering
    \begin{subfigure}[t]{0.48\textwidth}
        \centering
        \includegraphics[width=\textwidth]{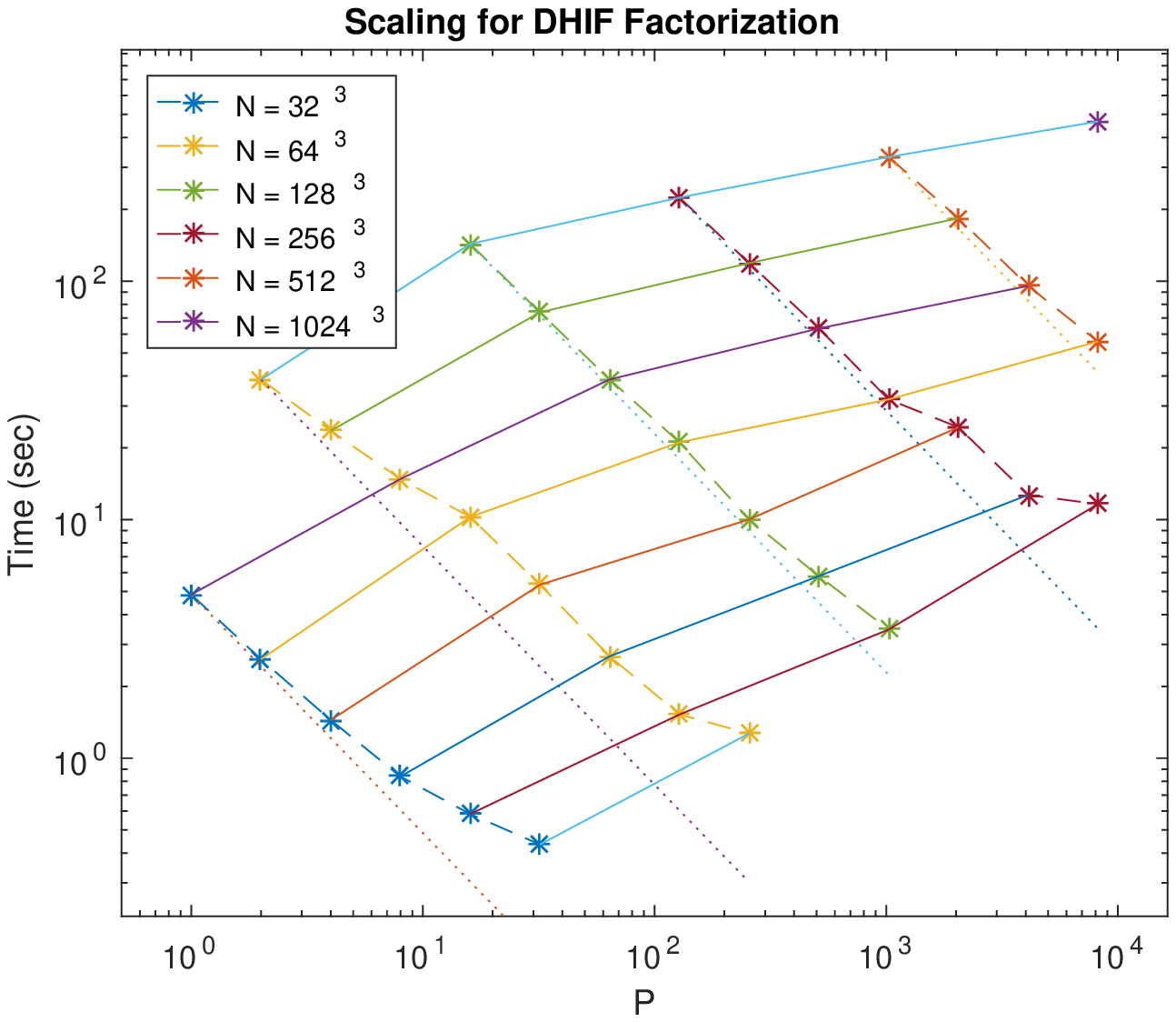}
        \caption{}
        \label{fig:ex1-facttss}
    \end{subfigure}
    ~
    \begin{subfigure}[t]{0.48\textwidth}
        \includegraphics[width=\textwidth]{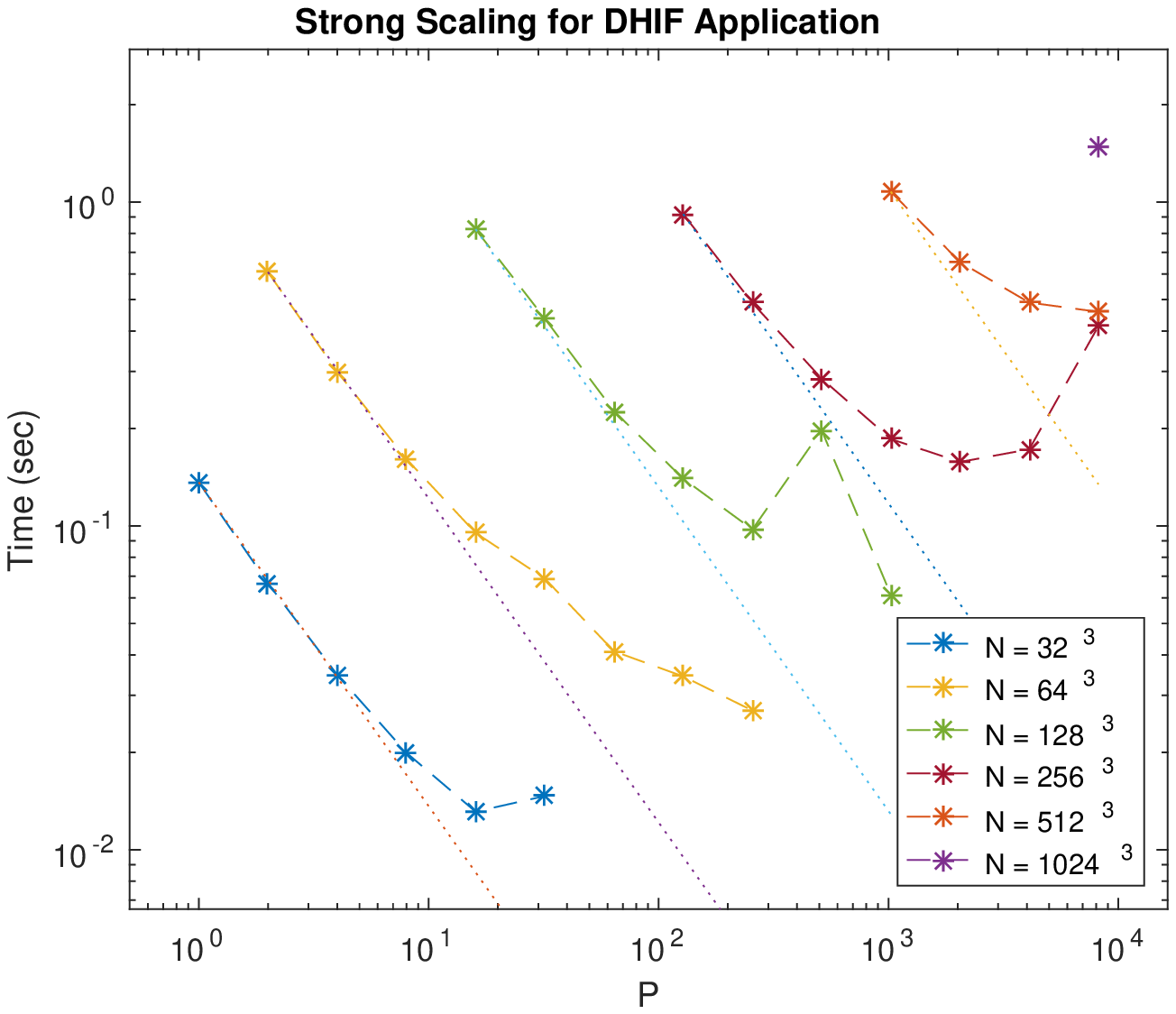}
        \caption{}
        \label{fig:ex1-apptss}
    \end{subfigure}

    \begin{subfigure}[t]{0.48\textwidth}
        \includegraphics[width=\textwidth]{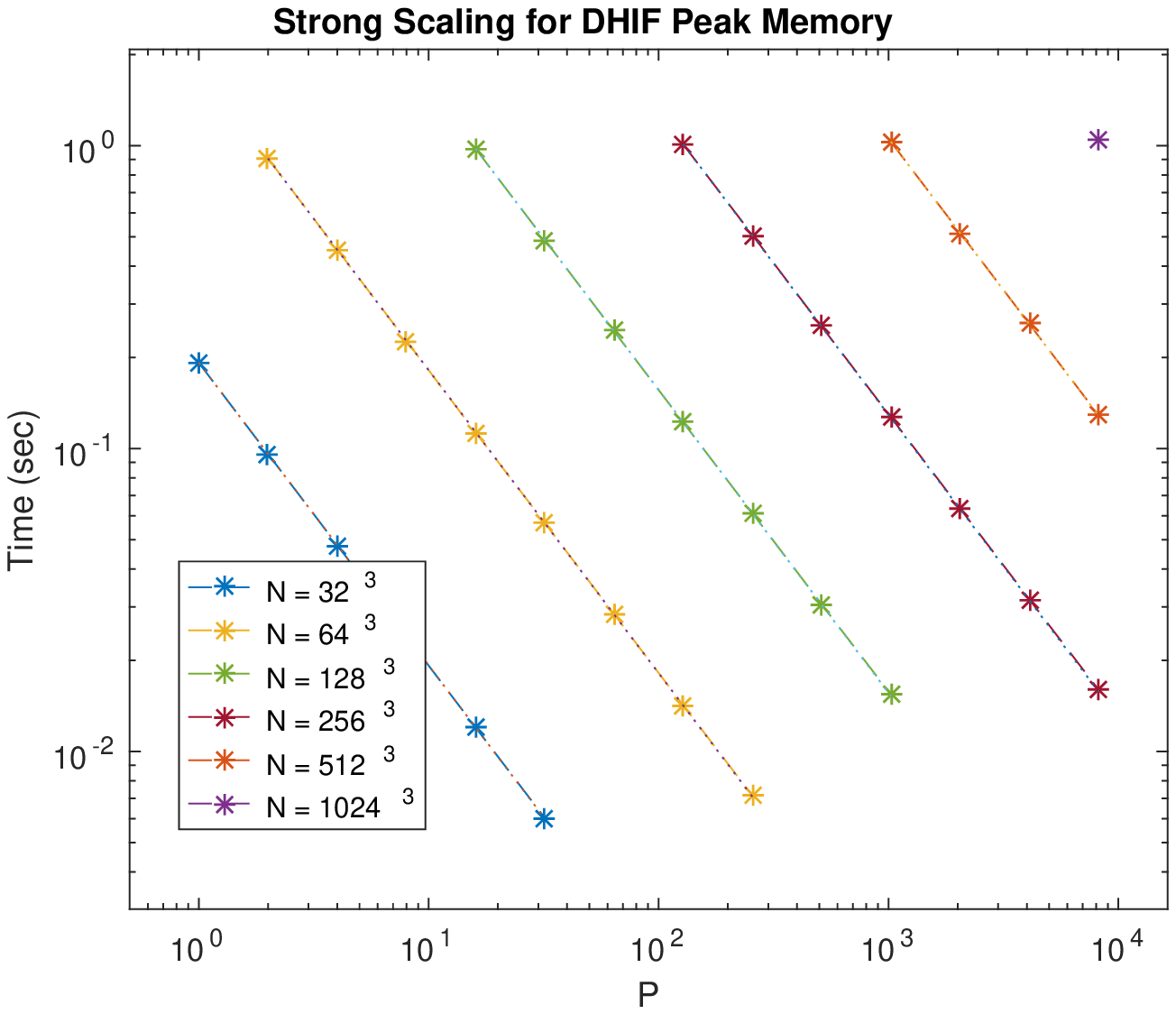}
        \caption{}
        \label{fig:ex1-memss}
    \end{subfigure}
    ~
    \begin{subfigure}[t]{0.48\textwidth}
        \includegraphics[width=\textwidth]{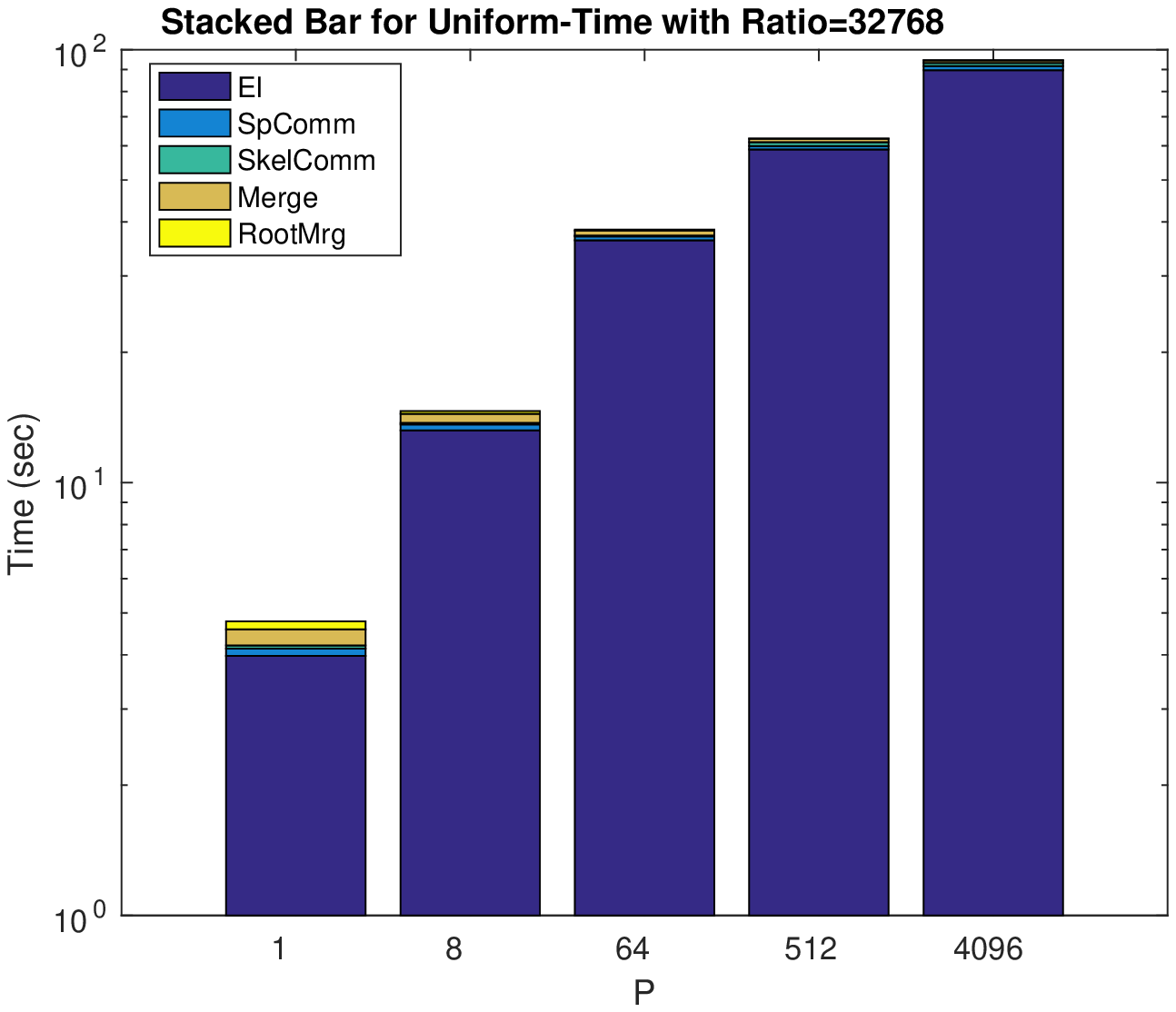}
        \caption{}
        \label{fig:ex1-time64sb}
    \end{subfigure}
    \caption{Example 1. (a) is the scaling plot for the DHIF
      factorization time, the solid lines indicate the weak scaling
      results, the dashed lines are the strong scaling results, and
      the dotted lines are the reference lines for perfect strong
      scaling, the line style apples to all figures below; (b) is the
      strong scaling for the DHIF application time; (c) is the strong
      scaling for the DHIF peak memory usage; (d) shows a stacked bar
      plot for factorization time for fixed ratio between the problem
      size and the number of processes.}
    \label{fig:ex1}
\end{figure*}

As shown in Table~\ref{tab:num-ex1}, given the tolerance $\epsilon =
10^{-3}$ the relative error remains well below this for all $N$ and
$P$. The number of skeleton points on the root level, $|\Sigma_L|$,
grows linearly as the one dimensional problem size increases.  The
empirical linear scaling of the root skeleton size strongly supports
the quasi linear scaling for the factorization, linear scaling for the
application, and linear scaling for memory cost.  The column labeled
with $m_f$ in Table~\ref{tab:num-ex1}, or alternatively
Figure~\ref{fig:ex1-memss}, illustrates the perfect strong scaling for
the memory cost.  Since the bottleneck for most parallel algorithms is
the memory cost, this point is especially important in practice.
Perfect distribution of the memory usage allows us to solve very large
problem on a massive number of processes, even through the
communication penalty on massive parallel computing would be
relatively large. The factorization time and application time show
good scaling up to thousands of processes.
Figure~\ref{fig:ex1-facttss} and Figure~\ref{fig:ex1-apptss} present
the strong scaling plot for the running time of factorization and
application respectively.  Together with
Figure~\ref{fig:ex1-time64sb}, which illustrates the timing for each
part of the factorization, we conclude that the communication cost
beside the parallel dense linear algebra (labeled with ``El'') remains
small comparing to the cost of the parallel dense linear algebra. It
is the parallel dense linear algebra part that stops the strong
scaling.  As it is also well know that parallel dense linear algebra
achieves good weak scaling, so does our DHIF implementation. Finally,
the last column of Table~\ref{tab:num-ex1} shows the number of
iterations for solving $Au=f$ using the GMRES algorithm with a
relative tolerance of $10^{-12}$ and with the DHIF as a
preconditioner. As the numbers in the entire column are equal to 6,
this shows that DHIF serves an excellent preconditioner with the
iteration number almost independent of the problem size.

{\bf \emph{Example 2.}}  The second example is a problem of
\eqref{eq:original-prob} with high-contrast random field $a(x)$ and
$b(x) \equiv 0.1$. The high-contrast random field $a(x)$ is defined as
follows,
\begin{enumerate}
\item Generate uniform random value $a_\vecj$ between $0$ and $1$
  for each discretization point;
\item Convolve the random value $a_\vecj$ with an isotropic
  three-dimensional Gaussian with standard deviation $1$;
\item Quantize the field via
  \begin{equation}
    a_\vecj = \left\{
    \begin{array}{ll}
      0.1, & a_\vecj \leq 0.5\\
      1000, & a_\vecj > 0.5\\
    \end{array}
    \right..
  \end{equation}
\end{enumerate}
The given tolerance is set to be $10^{-5}$.

\begin{figure}[htp]
  \centering
  \includegraphics[height=2in]{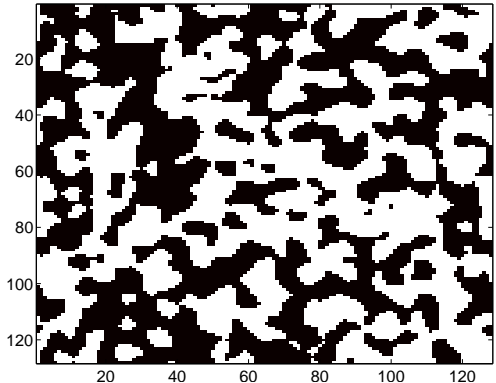}
  \caption{A slice in a random field realization of size $128^3$.}
  \label{fig:field}
\end{figure}

Figure~\ref{fig:field} shows a slice in a realization of the random
field.  The corresponding matrix $A$ is clearly of high-contrast.
Solving such a problem is harder than example 1 due to the raise of
the condition number.  The performance results of our algorithm are
presented in Table~\ref{tab:num-ex2}.  As we expect, the relative
error for solving is lower than that in Table~\ref{tab:num-ex1} and
the number of iterations in GMRES is higher.

\begin{table}
    \centering
    \begin{tabular}{crcrccrcrc}
\toprule
 $N$ & $P$ & $e_s$ & $|s_L|$ & $m_f$ & $t_f$ & $E^S_f$ & $t_s$ & $E^S_s$ & $n_{iter}$ \\
\toprule
\multirow{6}{*}{$32^3$}
& 1 & 3.02e-03 & 3865 & 2.00e-01 & 5.80e+00 & 100\% & 1.34e-01 & 100\% & 7 \\
& 2 & 3.39e-03 & 3632 & 9.31e-02 & 2.48e+00 & 117\% & 6.69e-02 & 100\% & 7 \\
& 4 & 2.69e-03 & 3934 & 5.13e-02 & 1.72e+00 &  84\% & 3.75e-02 &  89\% & 7 \\
& 8 & 3.18e-03 & 3660 & 2.37e-02 & 9.50e-01 &  76\% & 2.20e-02 &  76\% & 7 \\
& 16 & 3.13e-03 & 3693 & 1.24e-02 & 6.22e-01 &  58\% & 1.32e-02 &  63\% & 7 \\
& 32 & 3.00e-03 & 3744 & 6.42e-03 & 4.83e-01 &  38\% & 1.49e-02 &  28\% & 7 \\
\midrule
\multirow{8}{*}{$64^3$}
& 2 & 3.29e-03 & 8580 & 9.45e-01 & 4.33e+01 & 100\% & 6.15e-01 & 100\% & 7 \\
& 4 & 3.13e-03 & 8938 & 4.94e-01 & 2.91e+01 &  74\% & 3.10e-01 &  99\% & 7 \\
& 8 & 3.09e-03 & 9600 & 2.51e-01 & 1.98e+01 &  55\% & 1.68e-01 &  91\% & 7 \\
& 16 & 3.07e-03 & 8919 & 1.19e-01 & 1.27e+01 &  43\% & 9.86e-02 &  78\% & 7 \\
& 32 & 3.09e-03 & 9478 & 6.59e-02 & 6.99e+00 &  39\% & 7.89e-02 &  49\% & 7 \\
& 64 & 3.18e-03 & 9111 & 3.03e-02 & 3.17e+00 &  43\% & 4.90e-02 &  39\% & 7 \\
& 128 & 3.02e-03 & 9419 & 1.58e-02 & 2.15e+00 &  31\% & 3.31e-02 &  29\% & 7 \\
& 256 & 3.03e-03 & 9349 & 7.97e-03 & 1.60e+00 &  21\% & 3.66e-02 &  13\% & 7 \\
\midrule
\multirow{7}{*}{$128^3$}
& 16 & 3.16e-03 & 19855 & 1.07e+00 & 2.11e+02 & 100\% & 8.89e-01 & 100\% & 7 \\
& 32 & 3.09e-03 & 20487 & 5.58e-01 & 1.18e+02 &  90\% & 4.86e-01 &  91\% & 7 \\
& 64 & 3.06e-03 & 21345 & 2.78e-01 & 6.43e+01 &  82\% & 2.45e-01 &  91\% & 7 \\
& 128 & 3.10e-03 & 20344 & 1.37e-01 & 3.39e+01 &  78\% & 1.34e-01 &  83\% & 7 \\
& 256 & 3.07e-03 & 21152 & 7.43e-02 & 1.76e+01 &  75\% & 1.10e-01 &  51\% & 7 \\
& 512 & 3.07e-03 & 20779 & 3.51e-02 & 8.46e+00 &  78\% & 8.80e-02 &  32\% & 7 \\
& 1024 & 3.04e-03 & 21361 & 1.76e-02 & 5.38e+00 &  61\% & 6.31e-02 &  22\% & 7 \\
\midrule
\multirow{7}{*}{$256^3$}
& 128 & 3.11e-03 & 42420 & 1.14e+00 & 4.15e+02 & 100\% & 1.04e+00 & 100\% & 7 \\
& 256 & 3.12e-03 & 43828 & 5.91e-01 & 2.12e+02 &  98\% & 5.77e-01 &  90\% & 8 \\
& 512 & 3.11e-03 & 44126 & 2.90e-01 & 1.25e+02 &  83\% & 3.86e-01 &  67\% & 7 \\
& 1024 & 3.08e-03 & 43302 & 1.46e-01 & 6.31e+01 &  82\% & 2.12e-01 &  61\% & 7 \\
& 2048 & 3.09e-03 & 44131 & 7.78e-02 & 3.43e+01 &  76\% & 1.86e-01 &  35\% & 7 \\
& 4096 & 3.10e-03 & 43691 & 3.71e-02 & 1.96e+01 &  66\% & 2.28e-01 &  14\% & 7 \\
& 8192 & 3.10e-03 & 43952 & 1.85e-02 & 2.05e+01 &  32\% & 4.03e-01 &   4\% & 7 \\
\midrule
\multirow{4}{*}{$512^3$}
& 1024 & 3.11e-03 & 88070 & 1.16e+00 & 6.37e+02 & 100\% & 1.22e+00 & 100\% & 7 \\
& 2048 & 3.11e-03 & 88577 & 6.11e-01 & 3.47e+02 &  92\% & 6.84e-01 &  89\% & 8 \\
& 4096 & 3.11e-03 & 88757 & 3.03e-01 & 1.89e+02 &  84\% & 5.31e-01 &  58\% & 7 \\
& 8192 & 3.11e-03 & 85877 & 1.50e-01 & 1.02e+02 &  78\% & 6.20e-01 &  25\% & 7 \\
\midrule
\multirow{1}{*}{$1024^3$}
& 8192 & 3.11e-03 & 177323 & 1.18e+00 & 9.35e+02 & 100\% & 1.95e+00 & 100\% & 8 \\
\bottomrule
    \end{tabular}
    \caption{Example 2. Numerical results.}
    \label{tab:num-ex2}
\end{table}

\begin{figure*}[htp]
    \centering
    \begin{subfigure}[t]{0.48\textwidth}
        \centering
        \includegraphics[width=\textwidth]{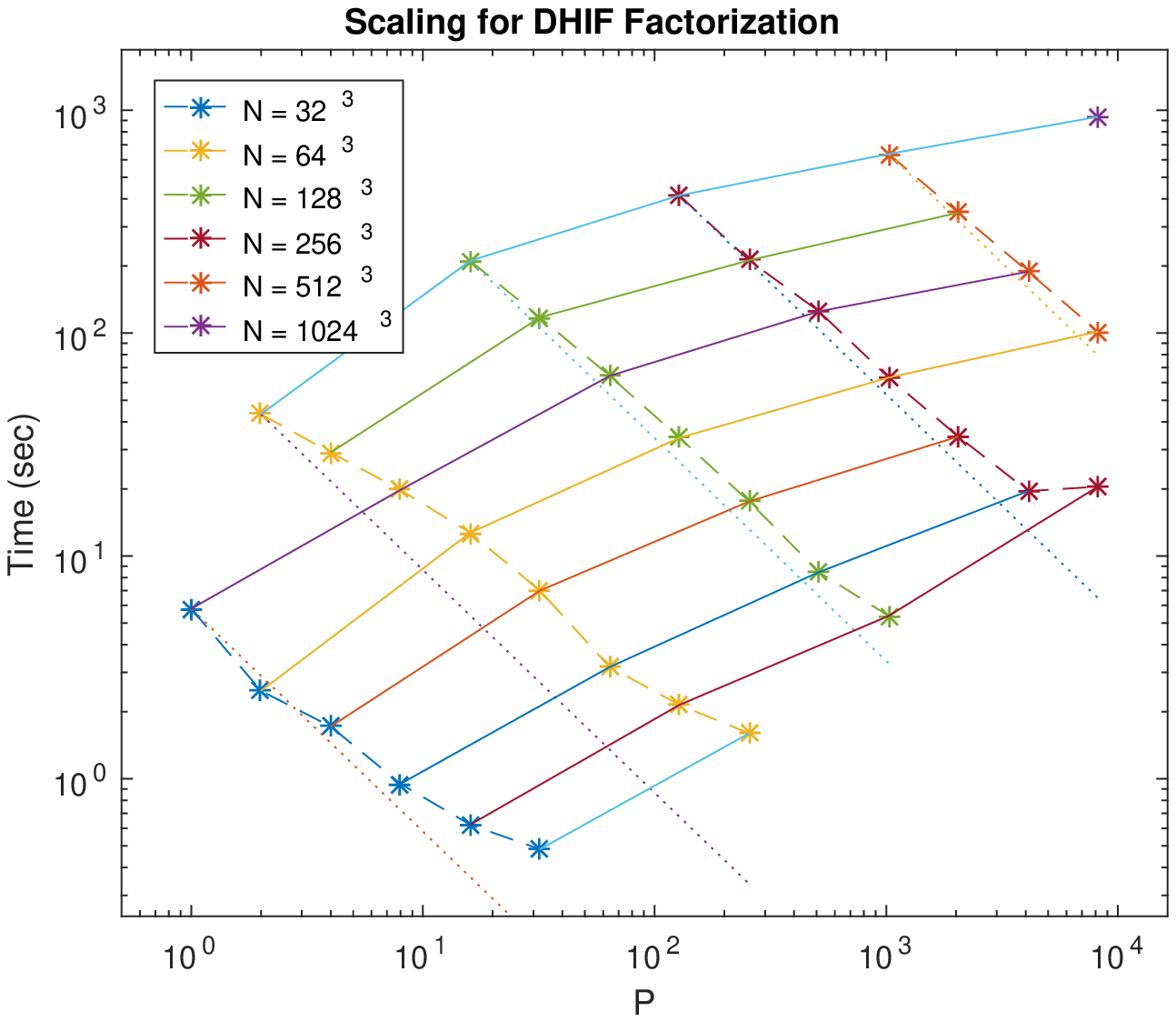}
        \caption{}
        \label{fig:ex2-facttss}
    \end{subfigure}
    ~
    \begin{subfigure}[t]{0.48\textwidth}
        \includegraphics[width=\textwidth]{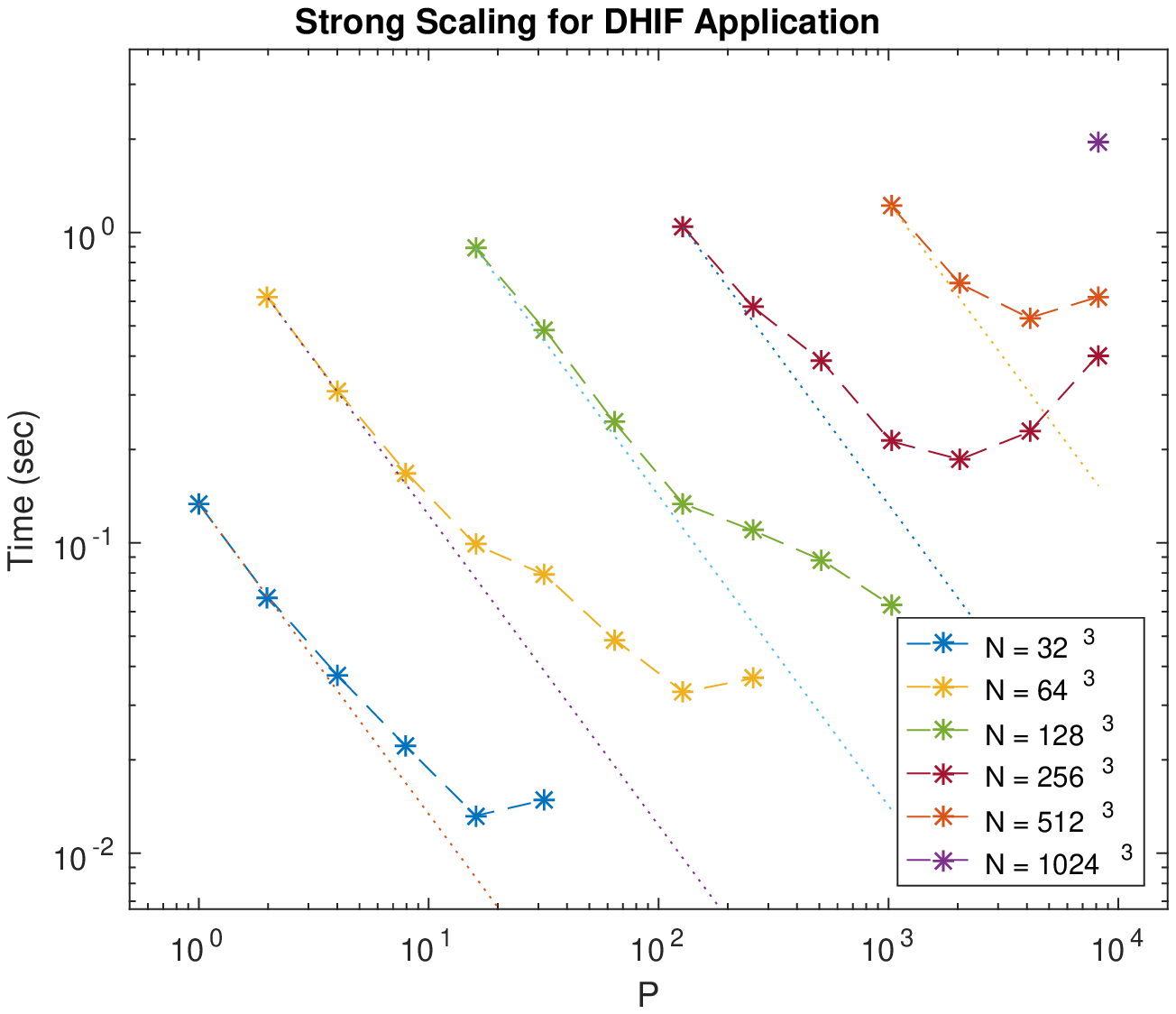}
        \caption{}
        \label{fig:ex2-apptss}
    \end{subfigure}

    \begin{subfigure}[t]{0.48\textwidth}
        \includegraphics[width=\textwidth]{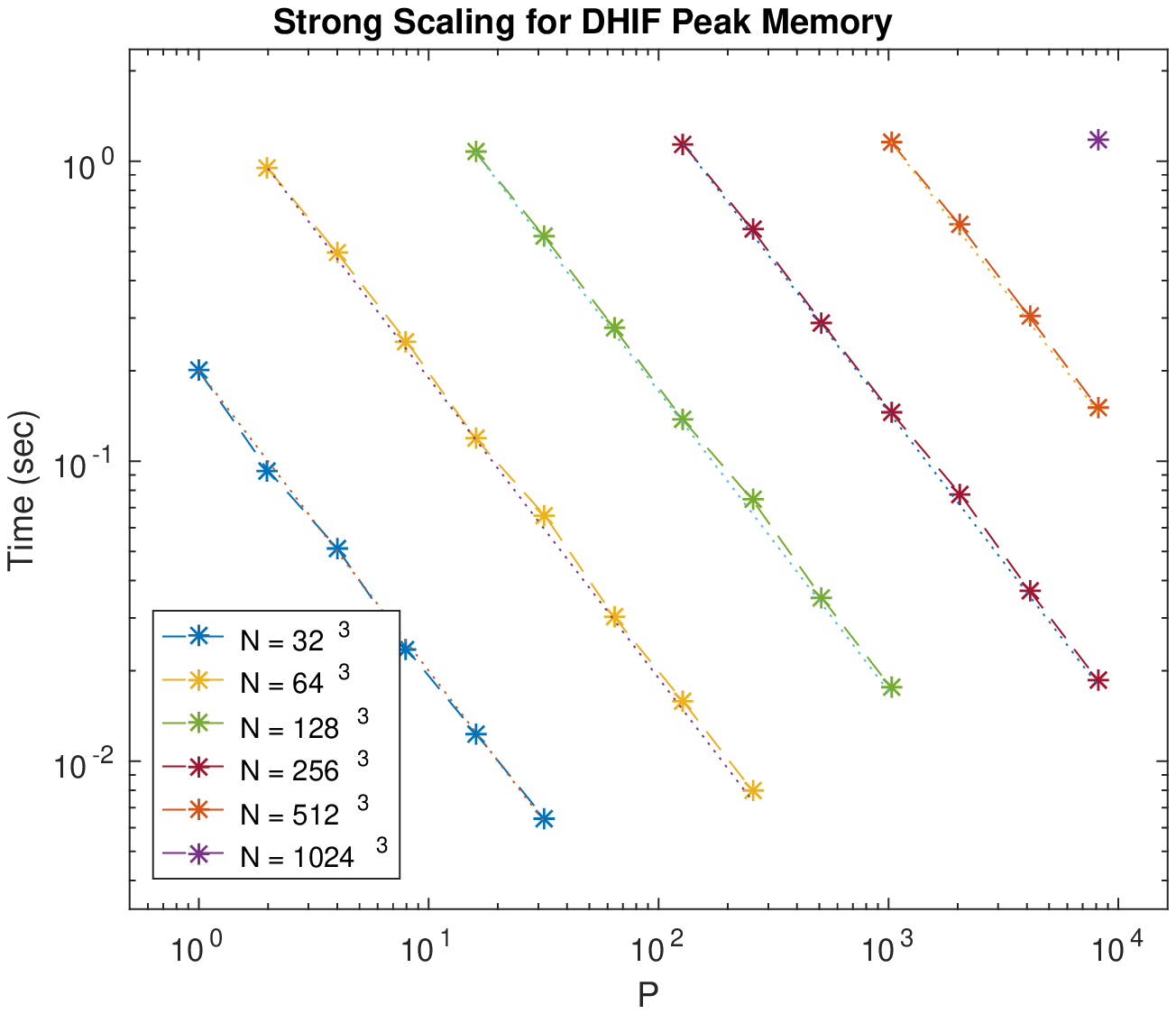}
        \caption{}
        \label{fig:ex2-memss}
    \end{subfigure}
    ~
    \begin{subfigure}[t]{0.48\textwidth}
        \includegraphics[width=\textwidth]{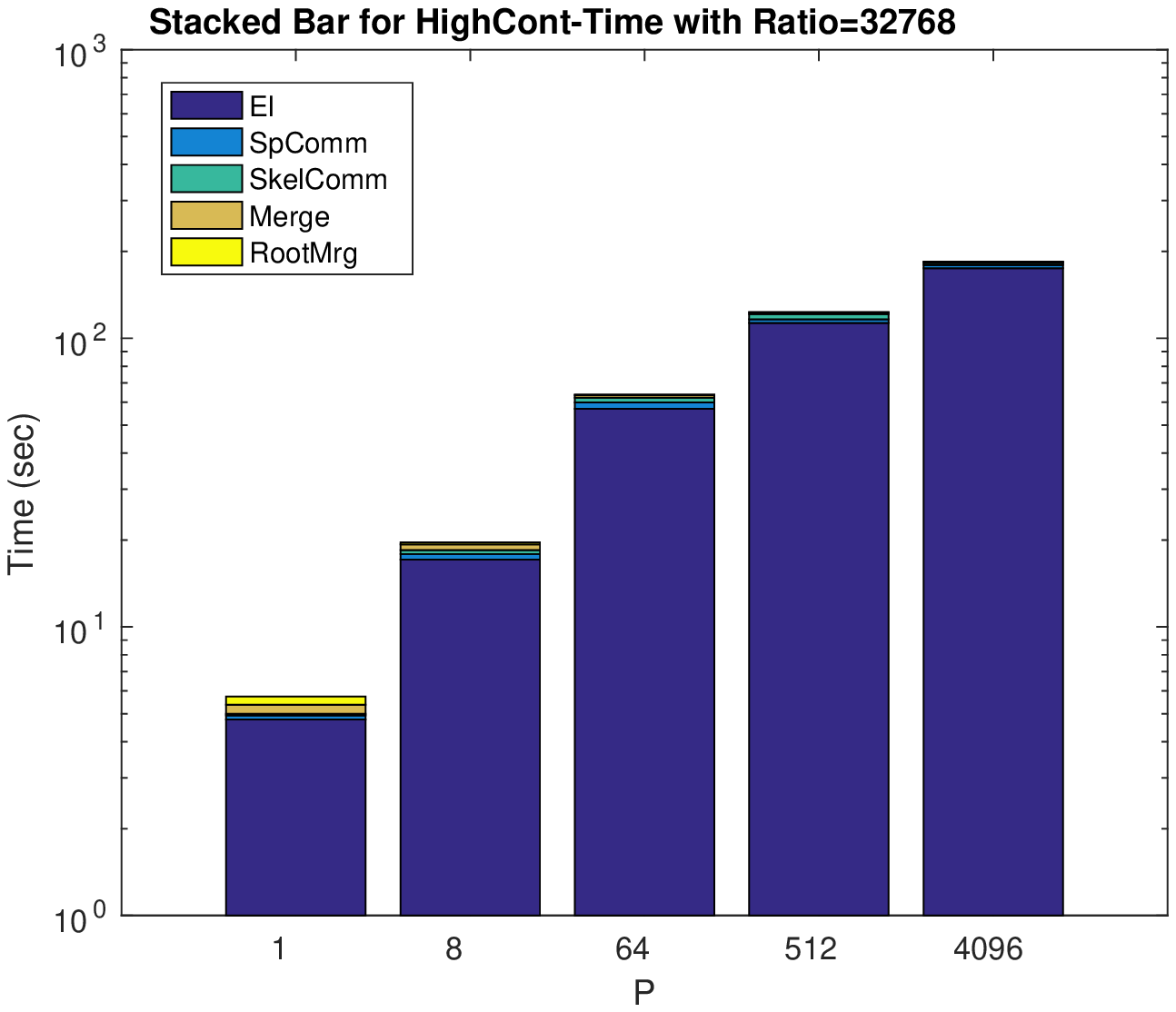}
        \caption{}
        \label{fig:ex2-time64sb}
    \end{subfigure}
    \caption{Example 2. (a) provides a scaling plot for
    DHIF factorization time;
    (b) is the strong scaling for DHIF application time;
    (c) is the strong scaling for DHIF peak memory usage;
    (d) shows a stacked bar plot for factorization time for fixed ratio
    between problem size and number of processes.}
    \label{fig:ex2}
\end{figure*}

Table~\ref{tab:num-ex2} and Figure~\ref{fig:ex2} demonstrate the
efficiency of the DHIF for high-contrast random field. Almost all
comments regarding the numerical results in Example 1 apply here.  To
focus on the difference between Example 1 and Example 2, the most
noticeable difference is about the relative error, $e_s$.  Though we
give a higher relative precision $\epsilon = 10^{-5}$, the relative
error for Example 2 is about $3\cdot 10^{-3}$, which is about ten
times larger than $e_s$ in Example 1. The reason for the decrease of
accuracy is most likely the increase of the condition number for
Example 2. This also increases the number of iterations in GMRES.
However, both $e_s$ and $n_{iter}$ remain roughly constant for varying
problem sizes. This means that DHIF still serves as a robust and
efficient solver and preconditioner for such problems. Another
difference is the number of skeleton points on the root level,
$|\Sigma_L|$.  Due to the fact that the field $a(x)$ is random, and
different rows in Table~\ref{tab:num-ex2} actually adopts different
realizations, the small fluctuation of $|\Sigma_L|$ for the same $N$
and different $P$ is expected. Overall $|\Sigma_L|$ still grows
linearly as $n = N^{1/3}$ increases. This again supports the
complexity analysis given above.

{\bf \emph{Example 3.}} The third example provides a concrete
comparison between the proposed DHIF and multigrid method
(hypre~\cite{Chow2006}). The problem behaves similar as example 2 without
randomness, \eqref{eq:original-prob} with high-contrast field $a(x)$ and
$b(x) \equiv 0.1$. The high-contrast field $a(x)$ is defined as follows,
\begin{equation}
    a(\vec{x}) = \left\{
        \begin{array}{ll}
            1000, & \sum_{i=1}^3 \lfloor \frac{x_i n}{7} \rfloor \equiv 0
            \pmod{2}\\
            0.1, & \sum_{i=1}^3 \lfloor \frac{x_i n}{7} \rfloor \equiv 1
            \pmod{2}\\
        \end{array}
        \right. ,
\end{equation}
where $n$ is the number of grid points on each dimension.

We adopt GMRES iterative method in both DHIF and hypre to solve the
elliptic problem to a relative error $10^{-12}$. The given tolerance
in DHIF is set to be $10^{-4}$. And SMG interface in hypre is used as
preconditioner for the problem on regular grids. The numerical results
for DHIF and hypre are given in Table~\ref{tab:num-ex3}.

\begin{table}[htb]
    \centering
    \begin{tabular}{cc|ccc|ccc}
        \toprule
        & & \multicolumn{3}{c|}{DHIF} & \multicolumn{3}{c}{hypre}\\
        $N$ & $P$ & $t_{setup} (sec)$ & $t_{solve} (sec)$ & $n_{iter}$
            & $t_{setup} (sec)$ & $t_{solve} (sec)$ & $n_{iter}$\\
        \toprule
        \multirow{2}{*}{$64^3$}
        &    8 &  15.27 & 18.10 & 21 &   0.29 &   9.67 &  67 \\
        &   64 &   2.46 &  3.45 & 21 &   1.47 &  17.37 &  60 \\
        \toprule
        \multirow{2}{*}{$128^3$}
        &   64 &  29.20 & 24.53 & 22  &  1.78 & 140.90 & 394 \\
        &  512 &   3.93 &  4.41 & 22  &  2.11 & 113.57 & 455 \\
        \toprule
        \multirow{2}{*}{$256^3$}
        &  512 &  59.66 & 26.33 & 21  &  4.11 & 258.22 & 492 \\
        & 4096 &  11.58 &  6.78 & 21  &  8.97 & 191.15 & 375 \\
        \bottomrule
    \end{tabular}
    \caption{Numerical results for DHIF and hypre. $t_{setup}$ is the
    setup time which is identical to $t_f$ in previous examples for DHIF,
    $t_{solve}$ is the total iterative solving time using GMRES, $n_{iter}$
    is the number of iterations in GMRES.} \label{tab:num-ex3}
\end{table}

As we can read from Table~\ref{tab:num-ex3}, there are a few advantages of
DHIF over hypre in the given settings. First, DHIF is faster than
hypre's SMG except for small problems with small numbers of processes.
And the number of iterations grows as the problem size grows in hypre, while
it remains almost the same in DHIF. In truely large problems, the advantages
of DHIF are more pronounced. Second, the scalability of DHIF appears to be
better than that of hypre's SMG. Finally, DHIF only requires powers of two
numbers of processes, whereas hypre's SMG requires powers of eight for 3D
problems.

\section{Conclusion}
\label{sec:Conclusion}

In this paper, we introduced the distributed-memory hierarchical
interpolative factorization (DHIF) for solving discretized elliptic
partial differential equations in 3D. The computational and memory
complexity for DHIF are
\begin{equation}
  \O\left(\frac{N\log N}{P}\right) \quad\text{and}\quad
  \O\left(\frac{N}{P}\right),
\end{equation}
respectively,
where $N$ is the total number of DOFs and $P$ is the number of processes.
The communication cost is
\begin{equation}
  \O\left(\sqrt{P}\log^3P\right)\alpha +
  \O\left(\frac{N^{2/3}}{\sqrt{P}}\right)\beta,
\end{equation}
where $\alpha$ is the latency, and $\beta$ is the inverse bandwidth.
Not only the factorization is efficient, the application can also be
done in $\O\left(\frac{N}{P}\right)$ operations.  Numerical examples
in Section~\ref{sec:Numerical Results} illustrate the efficiency and
parallel scaling of the algorithm. The results show that DHIF can be
used both as a direct solver and as an efficient preconditioner for
iterative solvers.

We have described the algorithm using the periodic boundary condition
in order to simplify the presentation. However, the implementation can
be extended in a straightforward way to problems with other type of
boundary conditions. The discretization adopted here is the standard
Cartesian grid. For more general discretizations such as finite
element methods on unstructured meshes, one can generalize the current
implementation by combining with the idea proposed in
\cite{Schmitz2014}.

Here we have only considered the parallelization of the HIF for
differential equations. As shown in \cite{Ho2015}, the HIF is also
applicable to solving integral equations with non-oscillatory
kernels. Parallelization of this algorithm is also of practical
importance.


{\bf Acknowledgments.}  Y. Li and L. Ying are partially supported by the
National Science Foundation under award DMS-1521830 and the U.S. Department
of Energy's Advanced Scientific Computing Research program under award
DE-FC02-13ER26134/DE-SC0009409. The authors would like to thank K. Ho,
V. Minden, A. Benson, and A. Damle for helpful discussions. We especially
thank J. Poulson for the parallel dense linear algebra package {\em
Elemental}. We acknowledge the Texas Advanced Computing Center (TACC)
at The University of Texas at Austin (URL: http://www.tacc.utexas.edu)
for providing HPC resources that have contributed to the research results
reported in the early versions of this paper. This research, in the
current version, used resources of the National Energy Research Scientific
Computing Center, a DOE Office of Science User Facility supported by
the Office of Science of the U.S. Department of Energy under Contract
No. DE-AC02-05CH11231.

\bibliographystyle{abbrv} \bibliography{ref}

\end{document}